\documentclass[review,onefignum,onetabnum]{siamart171218}

\usepackage{lipsum}
\usepackage{amsfonts}
\usepackage{graphicx}
\usepackage{epstopdf}
\ifpdf
  \DeclareGraphicsExtensions{.eps,.pdf,.png,.jpg}
\else
  \DeclareGraphicsExtensions{.eps}
\fi

\newsiamremark{remark}{Remark}
\newsiamremark{hypothesis}{Hypothesis}
\crefname{hypothesis}{Hypothesis}{Hypotheses}
\newsiamthm{claim}{Claim}

\headers{CEM-GMsFEM for sign-changing problems }{E. Chung, P. Ciarlet. X. Jin, and C. Ye}

\title{Multiscale Methods for wave propagation in materials with sign-changing coefficients\thanks{Submitted to the editors DATE.
\funding{The research of Eric Chung is partially supported by the Hong Kong RGC General Research Fund (Projects: 14305624 and 14305423 ).}}}

\author{Eric T. Chung\thanks{Department of Mathematics, The Chinese University of Hong Kong, Hong Kong
(\email{tschung@math.cuhk.edu.hk}).}
\and Patrick Ciarlet Jr.\thanks{POEMS, CNRS, INRIA, ENSTA Paris, Institut Polytechnique de Paris, 91120 Palaiseau, France
(\email{patrick.ciarlet@ensta.fr}).}
\and Xingguang Jin\thanks{Corresponding author. Department of Mathematics, The Chinese University of Hong Kong, Hong Kong
(\email{xgjin@math.cuhk.edu.hk}).}
\and Changqing Ye\thanks{Institute of Applied Physics and Computational Mathematics
  (\email{Ye\_Changqing@outlook.com}}).}

\usepackage{amsopn}

\makeatletter
\newcommand*{\addFileDependency}[1]{
  \typeout{(#1)}
  \@addtofilelist{#1}
  \IfFileExists{#1}{}{\typeout{No file #1.}}
}
\makeatother

\newcommand*{\myexternaldocument}[1]{%
    \externaldocument{#1}%
    \addFileDependency{#1.tex}%
    \addFileDependency{#1.aux}%
}

\ifpdf
\hypersetup{
  pdftitle={An Example Article},
  pdfauthor={D. Doe, P. T. Frank, and J. E. Smith}
}
\fi

\myexternaldocument{mms_supplement}
\usepackage{mathtools}

\allowdisplaybreaks

\usepackage{amssymb}
\usepackage[mathscr]{eucal}

\usepackage{mismath}
\nolinenumbers

\usepackage{bm} 

\usepackage{tikz}
\usetikzlibrary{math}
\usetikzlibrary{calc}

\usepackage{multirow}

\usepackage{makecell}
\setcellgapes{2pt}

\usepackage{booktabs}

\usepackage{siunitx}
\usepackage{algorithm}
\usepackage{algpseudocode}

\newcommand*\circled[1]{\tikz[baseline=(char.base)]{
    \node[shape=circle,draw,inner sep=1pt] (char) {\textbf{\small #1}};}}

\usepackage[skins]{tcolorbox}
\newtcolorbox{stepbox}[2][]{%
  enhanced,
  attach boxed title to top center={yshift=-3mm,yshifttext=-1mm},
  colframe=blue!75!black,
  colbacktitle=red!80!black,
  fonttitle=\bfseries,
  title=#2,#1
}

\newcommand{\Real}{\mathbb{R}}

\newcommand{\SubplotTag}[1]{\textbf{\small \textsf{#1}}}

\DeclareMathOperator{\Dist}{dist}
\DeclareMathOperator{\Span}{span}

\DeclareMathOperator{\Supp}{supp}

\DeclareMathOperator{\Diam}{diam}

\DeclarePairedDelimiter{\RoundBrackets}{(}{)}
\DeclarePairedDelimiter{\CurlyBrackets}{\{}{\}}

\newcommand{\dx}{\,\mathrm{d}}

\newtheorem{assumption}{Assumption}

\crefname{assumption}{assumption}{assumptions}
\Crefname{assumption}{Assumption}{Assumptions}
\DeclarePairedDelimiter{\Norm}{\lVert}{\rVert}
\crefname{equation}{}{}
\begin{document}
\maketitle
\begin{abstract}
From a mathematical perspective, the extraordinary properties of metamaterials are often reflected in the coefficients of the governing partial differential equations (PDEs). These coefficients may fall outside the assumptions of classical theory, particularly when the effective dielectric permittivity and/or magnetic permeability are negative. This situation can transform a coercive operator into a non-coercive one, potentially leading to ill-posedness. In this paper, we utilize the Constraint Energy Minimizing Generalized Multiscale Finite Element Method (CEM-GMsFEM), specifically designed for time-harmonic electromagnetic wave problems, where the construction of auxiliary spaces in the original CEM-GMsFEM is tailored to accommodate the sign-changing setting. Based on the framework of \texttt{T}-coercivity theory and resolution conditions, we establish the inf-sup stability and provide an a priori error estimate for the proposed method. The numerical results demonstrate the effectiveness and robustness of our approach in handling such sophisticated coefficient profiles.
\end{abstract}

\begin{keywords}
Negative index materials, T-coercivity, multiscale finite element method
\end{keywords}

\begin{AMS}
65M12, 65M15, 65N30
\end{AMS}
\section{Introduction}
Metamaterials have revolutionized the human understanding of materials science and engineering. Typically, metamaterials are engineered by arranging functional cells in a periodic structure, resulting in effective properties that exhibit behaviors markedly different from those of conventional materials. Veselago conceived the concept of negative refractive indices in 1967 \cite{Veselago1968}. Intuitively, light passing through these types of materials bends in the opposite direction compared to conventional materials. Much later, Pendry et al. proposed a viable design for experimental realization \cite{Pendry1999}, and Shelby, Smith, and Schultz \cite{Shelby2001} achieved the first successful fabrication in 2001. Another significant driving force behind the study of Negative-Index Metamaterials (NIM) was the potential for perfect lensing, which could surpass the diffraction limit of conventional lenses \cite{Pendry2000}. Given the considerable application potential associated with fabrication, in-silico simulations have become crucial for the design and optimization of NIMs.

From a mathematical perspective, the extraordinary properties of metamaterials are often reflected in the coefficients of the governing partial differential equations (PDEs), which may deviate from the assumptions of classical theory. Consider NIMs, where the electronic permittivity and magnetic permeability are both negative. When NIMs are embedded in a conventional medium, i.e. air, the governing PDEs exhibit coefficients that change sign globally. In a time-harmonic regime, this can transform a coercive operator into a non-coercive one, potentially leading to ill-posedness. This scenario poses significant challenges for both the theoretical analysis and numerical simulations. Over the past two decades, mathematicians have made considerable progress in understanding the rich insights of these PDEs, laying a robust foundation for numerical investigations. In particular, the $\texttt{T}$-coercivity theory introduced by Bonnet-Ben Dhia, Ciarlet, and Zwölf in their work \cite{Bonnet2016} is of particular importance, which establishes a unified framework for the well-posedness of PDEs with sign-changing coefficients.
The central idea of \texttt{T}-coercivity is to construct a bijective operator $\mathcal{T}$ such that the inf--sup condition is automatically satisfied. 
The construction of the operator $\mathcal{T}$ relies solely on the geometric configuration of the interfaces separating subdomains with coefficients of different signs. The \texttt{T}-coercivity theory guarantees the well-posedness of the model problem in situations where the Lax--Milgram theorem is not applicable. {In particular, well-posedness is ensured when the contrast ratio between the positive and negative coefficients is sufficiently large. The contrast ratio is defined as
$ \Upsilon = \frac{\sigma_{\min}^+}{\sigma_{\max}^-},$ where $\sigma_{\min}^+$ denotes the minimum positive coefficient and $\sigma_{\max}^-$ denotes the maximum absolute value of the negative coefficients.}  In \cite{BonnetBenDhia2012}, a sharp condition on the contrast ratio was established to guarantee well-posedness for certain classes of simple interface problems. The \texttt{T}-coercivity theory has since been extended to a broader range of problems, including Helmholtz-type equations \cite{CiarletJr.2012}, time-harmonic Maxwell systems \cite{BonnetBenDhia2014,Dhia2014}, eigenvalue problems \cite{Carvalho2017}, and mixed formulations \cite{Barre2023}.

The magnificent phenomena of NIMs are usually attributed to their interaction with common materials. Consequently, heterogeneity is a significant feature in metamaterial numerical simulations and should therefore be carefully addressed. Multiscale computational methods, which aim to capture the small-scale information of solutions due to heterogeneity at a reduced computational cost, are particularly well-suited for this task. Pioneered by Hou and Wu in \cite{Hou1997}, the methodology of incorporating model information into the construction of finite element spaces, coined as MsFEMs, has garnered significant attention. The evolution version, GMsFEMs, leverages spectral decomposition to perform dimension reduction for the online space, exhibiting superior performance when dealing with high-contrast and channel-like coefficient profiles \cite{Chung2014}. The first construction of multiscale bases capable of achieving the theoretically optimal approximation property for general rough coefficients is credited to M\aa{}lqvist and Peterseim in their celebrated work \cite{Maalqvist2014}. {This approach, known as the Localized Orthogonal Decomposition (LOD) method, employs quasi-interpolation operators to decompose the solution into macroscopic and microscopic components.} The combination of GMsFEMs and LOD led to the development of a CEM-GMsFEM by Chung, Efendiev, and Leung in \cite{Chung2018}, where “CEM” is the acronym for “Constraint Energy Minimizing”. Constrained Energy Minimizing Generalized Multiscale Finite Element Method (CEM-GMsFEM), which demonstrates superior performance in high-contrast problems, has been successfully
applied to various partial differential equations arising from practical applications \cite{Li2025,Zhao2023}. Constructing multiscale bases in CEM-GMsFEMs involves two steps: solving engineered generalized eigenvalue problems to obtain auxiliary spaces and localizing multiscale bases in a relaxed or constrained manner. Several primary successful results on Laplace-type PDEs with sign-changing coefficients have been achieved in \cite{ye2024} through the tailored construction of auxiliary spaces to accommodate the sign-changing setting. 

Despite these advancements, the current multiscale methods are still inadequate for NIM simulations, especially for electromagnetic wave problems. This claim is based on the following observed facts: (1) the current multiscale methods are primarily designed for Laplacian-type elliptic PDEs \cite{ye2024,Chaumont2021}. 
(2) Although with several numerical demonstrations, the high wave number regime remains challenging; (3) the sign-changing coefficients and the non-coercive nature undermine the numerical analysis. This paper aims to address the scalar electromagnetic wave equation in the time-frequency domain by developing GEM-GMsFEMs for NIM simulations.
Multiscale computational methods are initially designed to address heterogeneity in elliptic PDEs. However, accuracy can be compromised by high contrast coefficients in the PDEs. 
In the analysis section, we establish a resolution condition that aids in addressing the term involving the wave number. This condition allows us to absorb the wavenumber into the Laplacian term, facilitating the formulation of the inf-sup condition for both global and multiscale problems, thereby ensuring their well-posedness. The stability of GEM-GMsFEMs in the high contrast regime relies on auxiliary spaces, which are constructed by solving generalized eigenvalue problems. The well-posedness of CEM-GMsFEM is in fact guaranteed by utilizing the abstract framework of \texttt{T}-coercivity theory \cite{Dhia2014} to construct the inf-sup conditions. As emphasized earlier, CEM-GMsFEMs remain unexplored for conventional electromagnetic wave simulations. These observations motivate us to develop a novel multiscale method, starting with common electromagnetic wave simulations and progressing to NIM simulations.

To the best of our knowledge, efforts to apply multiscale computational methods to electromagnetic wave problems are currently scarce. An exception can be found in the work by Chung and Ciarlet in \cite{Chung2013}, where they employed the staggered discontinuous Galerkin method \cite{Chung2013a,Chung2012}.
Note that contrast ratios play a significant role in the \texttt{T}-coercivity theory to assess the well-posedness of the model problem.
In comparison, the proposed method presented in this paper is more general, as it can handle a broader range of coefficient profiles, including those with high contrast ratios.

This paper is organized as follows. In \cref{sec:preliminaries}, we introduce the model problem and present the \texttt{T}-coercivity theory. The construction of multiscale bases in the proposed method is detailed in \cref{sec:CEM-GMsFEM}.
To validate the performance of the proposed method, \cref{sec:Numerical} presents numerical experiments conducted on four different models.
All theoretical analysis for the proposed method is gathered in \cref{sec:anal}.
Finally, in \cref{sec:conclusions}, we conclude the paper.

\section{Preliminaries} \label{sec:preliminaries}              
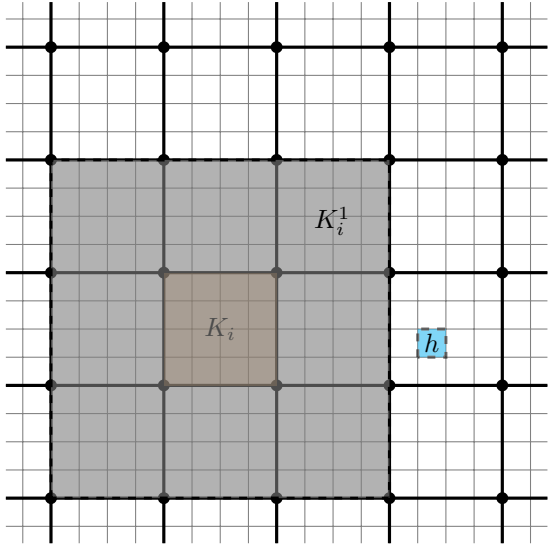
\begin{figure}[!ht]
\centering
\begin{tikzpicture}[scale=1.5]
\draw[step=0.25, gray, thin] (-0.4, -0.4) grid (4.4, 4.4);
\draw[step=1.0, black, very thick] (-0.4, -0.4) grid (4.4, 4.4);
\foreach \x in {0,...,4}
\foreach \y in {0,...,4}{
\fill (1.0 * \x, 1.0 * \y) circle (1.5pt);
}
\fill[brown, opacity=0.4] (1.0, 1.0) rectangle (2.0, 2.0);
\node at (1.5, 1.5) {$K_i$};
\draw [dashed, very thick, fill=gray, opacity=0.6] (0, 0) rectangle (3, 3);
\node[above right] at (2.25, 2.25) {$K_i^1$};
\draw [dashed, very thick, fill=cyan, opacity=0.5] (3.25, 1.25) rectangle (3.5, 1.5);
\node at (3.375, 1.375) {$h$};
\end{tikzpicture}
\caption{An illustration of the two-scale mesh, a fine element $h$, a coarse element $K_i$ and its oversampling coarse element $K_i^m$ with the oversampling layer $m=1$.}
\label{fig:grid}
\end{figure}

Let $\Omega \subset \mathbb{R}^2$ be a bounded open domain, decomposed into two non-overlapping subdomains $\Omega^+$ and $\Omega^-$ with {$\Gamma$ as the interfaces.} Here, we adopt the notations for all quantities $\sigma,\sigma^{-1},c,c^{-1}\in L^{\infty}(\Omega)$ defined on $\Omega$ such that 
\begin{equation*}
\begin{cases}
x\in\Omega^+,\sigma(x)>0:\sigma_{\mathup{max}}^+=\sup\sigma(x),\sigma_{\mathup{min}}^+=\inf\sigma(x).\\
x\in\Omega^-,\sigma(x)<0:\sigma_{\mathup{max}}^-=\sup\abs{\sigma(x)},\sigma_{\mathup{min}}^-=\inf\abs{\sigma(x)}.\
\end{cases}
\end{equation*}
We consider the following Helmholtz-like equation with sign-changing coefficients $\sigma$ and $c$ in the bounded space domain $\Omega\subset\Real^2$ imposed by homogeneous boundary conditions:
\begin{equation}\label{eq:ell1}
-\nabla\cdot(\sigma\nabla u)-k^2 cu=f,
\end{equation}
where $k\in\Real$ is a positive wavenumber, $ f\in L^2(\Omega)$ represents a harmonic source. Redefine the Sobolev space $V=H_0^1(\Omega)$ 
and write the boundary value problem \cref{eq:ell1} in a variational form: find a solution $u\in V$ such that for all $v\in V$,
\begin{equation}\label{eq:weak}
\begin{aligned}
\int_\Omega \sigma\nabla u\cdot\nabla v \di x-k^2\int_{\Omega}c uv\di x &=\int_{\Omega} fv \di x, &&\forall v\in {V},\\
\end{aligned}
\end{equation}
To simplify the notations, the sesquilinear form $\mathcal{B}: V\times V\rightarrow\mathbb{R}$ satisfies
 \begin{equation}\label{bform}
\mathcal{B}(u,v) :=(\sigma\nabla u, \nabla v)-k^2(c u,v),
 \end{equation}
where $\RoundBrackets*{u,v}=\int_{\Omega}uv\di x$. Then we rewrite \cref{eq:weak} in the following
\begin{equation}\label{vf}
\mathcal{B}(u,v)=(f,v)=l(v).
\end{equation}
It is convenient
to introduce a notation for a norm as
$\| v \|_{\tilde{a}, \omega}
\coloneqq
\left( \int_{\omega} |\sigma| \, |\nabla v|^2 \, \mathrm{d}x \right)^{1/2},
$
where $\omega \subset \Omega$ is a subdomain. 
{Moreover, when $\omega = \Omega$, we simply write 
$\| \cdot \|_{\tilde{a}}$ instead of $\| \cdot \|_{\tilde{a}, \Omega}$.} Based on the definition of \texttt{T}-coercivity: If there exists a bijective operator $\mathcal{T}:V\rightarrow V$, such that $\forall v\in V$, $\exists \gamma>0$, a bilinear form $a$ satisfies $|a(v, \mathcal{T}v)|\geq\gamma\Norm{v}_{\tilde{a}}^2$.
If we split form $\mathcal{B}$ as $\mathcal{B}=a+b$ where $a=(\sigma\nabla u, \nabla v)$ and $b=-k^2(cu,v)$, then the \texttt{T}-coercivity of the bilinear form $a(v, \mathcal{T}v)$ can be obtained under the sufficiently large high-contrast ratio $\Upsilon=\frac{\sigma_{\text{min}}^+}{\sigma_{\text{max}}^-}$ in \cite{ye2024}, and the construction of $\mathcal{T}$ {is} given by
\begin{equation}\label{eq:T operator}
 \mathcal{T}v=\begin{cases}
    v_1,                  & \text{ in } \Omega^+, \\
    -v_2+2\mathcal{R}v_1, & \text{ in } \Omega^-,
  \end{cases}
\end{equation}
where
\begin{equation}\label{eq:split of v}
  v=\begin{cases}
    v_1, & \text{ in } \Omega^+, \\
    v_2, & \text{ in } \Omega^-.
  \end{cases}
\end{equation}
 Here, we assume the existence of a bounded map $\mathcal{R}\colon V^+\rightarrow V^-$ with $\mathcal{R}v|_{\Gamma}=v|_{\Gamma}$ in the sense of Sobolev traces, and the second term $b(\cdot, \mathcal{T}\cdot)$ is a compact perturbation due to the compacting embedding theorem \cite{Bonnet2016}, then the weak form \cref{vf} is well-posed if and only if $l(v)=0$ implies $u=0$.
 {
\begin{remark}\label{remark1}
Note that the derivation in \cref{eq:T operator} involves switching the
``positive'' subspace to the ``negative'' one, according to the definition
of $\mathcal{T}\colon V^+ \rightarrow V^-$. 
Alternatively, one may reverse this mapping from ``negative'' to ``positive''. 
In this case, the \texttt{T}-coercivity of the bilinear form $a$ can still be established,
provided that the ratio $\sigma_{\mathup{min}}^- / \sigma_{\mathup{max}}^+$ 
is sufficiently large (or $\sigma_{\mathup{max}}^+/\sigma_{\mathup{min}}^-$ is sufficiently small), as discussed in \cite{BonnetBenDhia2012}.
\end{remark}}
Based on \cite{Chesnel2013}, and as shown in \cref{fig:grid}, let $\mathcal{K}_H$ be a standard quadrilateral partition of the domain $\Omega$ with mesh size $H$. We refer to this partition as coarse grids, and the produced elements as the coarse elements. Each coarse element $K\in\mathcal{K}_H$ is further partitioned into a union of connected fine grid blocks. We denote the fine-grid partition as $\mathcal{K}_h$ with $h>0$ being the fine grid size. For each $K_i\in \mathcal{K}_H$ with $1\leq i\leq N$, $N$ is the number of coarse elements. we define an oversampled domain $K_i^m(m\geq 1)$ in the following 
$$
K_i^m: =\text{int}\CurlyBrackets*{\bigcup_{K\in\mathcal{K}_H}\CurlyBrackets*{  \overline{K_i^{m-1}}\bigcap\overline{K}\neq \varnothing}\bigcup\overline{K_i^{m-1}}},
$$
where $\text{int}(S)$ and $\overline{S}$ represent the interior and the closure of a set $S$, and the initial value $K_i^0=K_i$ for each element.
\section{The construction of the CEM-GMsFEM basis functions}\label{sec:CEM-GMsFEM}
To facilitate the construction of the eigenvalue problems and the subsequent analysis, we can assume that the value of 
$c$ in the equation \cref{eq:ell1} is proportional to $\sigma$. In the original CEM-GMsFEM, a generalized eigenvalue problem,
\begin{equation}\label{eq:old aux}
\begin{aligned}
&\text{find } \lambda \in \mathbb{R} \text{ and } v \in H^1(K_i)\setminus\CurlyBrackets{0} \\
&\text{ s.t. }\forall w\in H^1(K_i),
\int_{K_i} \sigma \nabla v \cdot \nabla w \di x = \lambda \int_{K_i} \mu_\mathup{msh} \text{diam}(K_i)^{-2} cv w \di x,
\end{aligned}
\end{equation}
is solved on each coarse element $K_i$, where $\mu_\mathup{msh}$ is a generic positive constant that depends on the mesh quality.
Then the local auxiliary space $V_i^\mathup{aux}$ is formed by collecting leading eigenvectors.
However, recalling that $\sigma$ is not uniformly positive, the left-hand bilinear form in \cref{eq:old aux} is not positive semidefinite, and similarly, the right-hand bilinear form in \cref{eq:old aux} is not positive definite.
Then, determining leading eigenvectors via \cref{eq:old aux} is problematic since the eigenvalues could be negative.
To address this issue, we instead construct the following generalized eigenvalue problem on each $K_i$, which forms the first step of the proposed method:
\begin{equation}\label{eq:new aux}
\begin{aligned}
 & \text{find } \lambda \in \mathbb{R} \text{ and } v \in H^1(K_i)\setminus\CurlyBrackets{0} \\
 &\text{ s.t. } \forall w \in H^1(K_i), 
  \int_{K_i} \abs{\sigma} \nabla v \cdot \nabla w \di x = \lambda \int_{K_i} \mu_\mathup{msh} \text{diam}(K_i)^{-2}\abs{c}v w \di x.
\end{aligned}
\end{equation}
Upon solving the $l_*$ leading eigenvectors in \cref{eq:new aux}, denoted as $\psi_i^j$ with $1 \leq j \leq l_*$, we can construct the local auxiliary space $V_i^\mathup{aux} \subset L^2(K_i)$ as $\Span\{\psi_i^j\}$ for $1 \leq i \leq N$.
The global auxiliary space $V^\mathup{aux} \subset L^2(\Omega)$ is defined as $V^\mathup{aux}=\oplus_{i=1}^{N}V_i^\mathup{aux}$, where $V_i^\mathup{aux}\subset L^2(\Omega)$ is the space by performing zero-extension of functions in $V_i^\mathup{aux}$. We introduce several notations for future reference. Let $\mu$ be a function in $L^\infty(\Omega)$ satisfying
\[
  \mu|_{K_i}=\mu_\mathup{msh}\text{diam}(K_i)^{-2}c|_{K_i}
\]
for all $K_i \in \mathcal{K}_H$.
Next, for a subdomain $\omega$ in $\Omega$, we define two bilinear forms:
\[
a(v, w)_{\omega} \coloneqq \int_{\omega} \sigma \nabla v \cdot \nabla w \di x \text{ and } s(v, w)_{\omega} \coloneqq \int_{\omega} \mu v w \di x.
\]
Similarly to the definition of $\norm{v}_{\tilde{a},\omega }$, we define $\norm{v}_{\tilde{s},\omega }$ as $(\int_{\omega} \abs{\mu} \abs{v}^2 \di x)^{1/2}$.
Again, we drop the subscript $\omega$ if $\omega=\Omega$.
Additionally, we define the orthogonal projection operator $\pi\colon L^2(\Omega) \rightarrow L^2(\Omega)$ under the norm $\norm{\cdot}{\tilde{s}}$, with $V^\mathup{aux}=\im \pi$. The following Lemma \ref{inter} demonstrates the properties of the global projection $\pi$, which will be frequently utilized in the analysis. Its proof is a straightforward application of the properties of eigenspace expansions.
We also use the shorthand notation $V_i^m$ as $H^1_0(K_i^m)$, and we always implicitly identify $V_i^m$ as a subspace of $V$. The functions in \( V^\text{aux} \) may not be continuous in \( \Omega \), and thus \( V^\text{aux} \) cannot be used as a finite element space. In the next section, \cref{mbf}, we will construct a multiscale basis \( \phi_{i,m}^j \) in \( V_i^m \), corresponding to \( \psi_{i}^j \) in \( V_i^\text{aux} \), which is obtained in equation \cref{eq:new aux}. Before that, we will list a collection of estimates related to the \( \mathcal{T} \) operator, which will help pave the way for proving the well-posedness of the global multiscale basis function in equation \cref{msbasis}.
\begin{lemma}\label{inter}
For each $K_i \in \mathcal{K}_H$ and for all $v \in H^1(K_i)$, the following estimates hold:
\begin{equation}\label{pi1}
\| v - \pi_i v \|_{\tilde{s}, K_i}^2
\leq \frac{\| v \|_{\tilde{a}, K_i}^2}{\lambda_i^{l_*+1}}
\leq \Lambda^{-1} \| v \|_{\tilde{a}, K_i}^2,
\end{equation}
where $\Lambda = \min_{1 \leq i \leq N} \lambda_i^{l_*+1},$ and {$\lambda_i^{l_*+1}$ denotes the $(l_*+1)$-th eigenvalue of the generalized eigenvalue problem \cref{eq:new aux}}. Moreover,
\begin{equation}\label{pi2}
\| \pi_i v \|_{\tilde{s}, K_i}^2
= \| v \|_{\tilde{s}, K_i}^2
- \| v - \pi_i v \|_{\tilde{s}, K_i}^2
\leq \| v \|_{\tilde{s}, K_i}^2.
\end{equation}
\end{lemma}
\begin{lemma}\label{well}
For any $v \in V$, where $\mathcal{T}$ and $\mathcal{R}$ are the operators defined in \cref{eq:T operator}, the following estimates hold:
\begin{align}
a(v, \mathcal{T}v)
& \geq
\RoundBrackets*{1 - \|\mathcal{R}\|_1 / \sqrt{\Upsilon}}
\|v\|_{\tilde{a}}^2,
\label{eq:coer a}
\\
s(\pi v, \pi \mathcal{T}v)
& \geq
\RoundBrackets*{1 - \|\mathcal{R}\|_0 / \sqrt{\Upsilon}}
\|\pi v\|_{\tilde{s}}^2
-
\Lambda^{-1}
\|\mathcal{R}\|_0/\sqrt{\Upsilon}
\|v\|_{\tilde{a}}^2,
\label{eq:coer s}
\\
\|\mathcal{T}v\|_{\tilde{a}}
& \leq
\max\CurlyBrackets*{
\RoundBrackets*{1 + 8\|\mathcal{R}\|_1^2 / \Upsilon}^{1/2},
\sqrt{2}
}
\|v\|_{\tilde{a}},
\label{eq:bound a}
\\
\|\mathcal{T}v\|_{\tilde{s}}
& \leq
\max\CurlyBrackets*{
\RoundBrackets*{1 + 8\|\mathcal{R}\|_0^2 / \Upsilon}^{1/2},
\sqrt{2}
}
\|v\|_{\tilde{s}}.
\label{eq:bound s}
\end{align}
Here, for any $w \in V(\Omega^+)$,
\begin{equation}\label{eq:norm1 of R}
\|\nabla \mathcal{R}w\|_{0,\Omega^-}
\leq
\|\mathcal{R}\|_1 \|\nabla w\|_{0,\Omega^+},
\qquad
\|\mathcal{R}w\|_{0,\Omega^-}
\leq
\|\mathcal{R}\|_0 \|w\|_{0,\Omega^+},
\end{equation}
{where $\|\mathcal{R}\|_0$ and $\|\mathcal{R}\|_1$ are positive constants, and 
$\|v\|_{0,\omega} := \left( \int_{\omega} |v|^2 \, \mathrm{d}x \right)^{1/2}$
for any subdomain $\omega \subset \Omega$.}
\end{lemma}
\begin{proof}
{The detailed proof of this lemma can be found in \cite{ye2024}.} 
To facilitate the analysis, we introduce the following notations for constants that are from \cref{eq:coer a,eq:coer s,eq:bound a,eq:bound s}:
\begin{align*}
  C_0 & \coloneqq \max\CurlyBrackets*{\RoundBrackets*{1+8\norm{\mathcal{R}}_1^2 / \Upsilon}^{1/2},\ \sqrt{2}},                                                       \\
  C_1 & \coloneqq \max\CurlyBrackets*{\RoundBrackets*{1+8\norm{\mathcal{R}}_0^2 / \Upsilon}^{1/2},\ \sqrt{2}},                                                       \\
  C_2 & \coloneqq \min\CurlyBrackets*{1-\norm{\mathcal{R}}_0/\sqrt{\Upsilon}- \Lambda^{-1}\norm{\mathcal{R}}_0/\sqrt{\Upsilon},\ 1-\norm{\mathcal{R}}_1/\sqrt{\Upsilon}}.
\end{align*}
We can see that $C_0$, $C_1$, and $C_2$ can all be bounded from above and below, provided that $\Upsilon$ is sufficiently large and $\Lambda^{-1}$ is sufficiently small.
\end{proof}
\subsection{Multiscale basis functions}\label{mbf}
For each coarse element $K_i\in\mathcal{K}_H$ and its oversampling domain $K_i^m\subset\Omega$ by enlarging $K_i$ for $m$ coarse grid layers, we define the multiscale basis function $\phi_{i,m}^j\in V_0(K_i^m)$. Find $\phi_{i,m}^j\in V_0(K_i^m)$ such that
\begin{equation}\label{msbasis1}
\mathcal{B}(\phi_{i,m}^j,w)_{K_i^m}+s(\pi \phi_{i,m}^j,\pi w)_{K_i^m}=s(\psi_{i}^j,\pi w)_{K_i},\quad\forall w\in V_0(K_i^m),
\end{equation}
where $V_0(K_i^m)$ is the subspace of $V(K_i^m)$ with zero trace on $\partial K_i^m$ and $V(K_i^m)$ is the restriction of $V$ in $K_i^m.$ 
The global multiscale basis function $\phi_i^j\in V$ is defined in a similar way, 
\begin{equation}\label{msbasis}
\mathcal{B}(\phi_{i}^j,w)+s(\pi \phi_{i}^j,\pi w)=s(\psi_{i}^j,\pi w)_{K_i},\quad\forall w\in V(K_i).
\end{equation}
{Note that the variational problems defined in \cref{msbasis1} and \cref{msbasis} may become ill-defined due to the sign-changing coefficients. 
However, in practice, they are solved on the fine mesh $\mathcal{K}_h$. 
To overcome this difficulty, we introduce modified formulations of both the local and global problems, inspired by \cref{msbasis1} and \cref{msbasis} but adapted to the $\texttt{T}$-coercivity setting. The well-posedness of these revised variational problems will be rigorously established in the subsequent analysis 
\cref{sec:anal}. 
}
{We first define the associated global operator $\mathcal{G}^\infty_i\colon L^2(\Omega) \rightarrow V=H^1_0(\Omega)$ corresponding to the coarse element $K_i$ via the following variational problem:
\begin{equation}\label{eq:global operator}
  \text{Find } \mathcal{G}_i^\infty \psi \in V \text{ s.t. } \forall w \in V,\
  \mathcal{B}(\mathcal{G}_i^\infty \psi, w) + s(\pi\mathcal{G}_i^\infty \psi, \pi w) = s(\pi\psi, \pi w)_{K_i}.
\end{equation}}
By defining $\mathcal{G}^{\infty} \psi = \sum_{i=1}^{N} \mathcal{G}_i^{\infty} \psi$ and $V_{\text{glo}} = \mathrm{Im}(\mathcal{G}^{\infty})$, the global solution $u_{\text{glo}} = \mathcal{G}^{\infty} \psi$ is then defined as the solution to the following variational problem:
\begin{equation}\label{eq:global solution}
\text{find } u_{\mathup{glo}}\in V_{\mathup{glo}},\ \text{ s.t. }\forall v\in V_{\mathup{glo}},\ \mathcal{B}(u_{\mathup{glo}},v)=(f,v).
\end{equation}
\begin{lemma}\label{lm3.2}
For any $v_{\mathup{glo}}\in V_{\mathup{glo}}$ and $w\in \tilde{V}$, it holds $\mathcal{B}(v_{\mathup{glo}},w)=0.$ If $w\in V$ and $\mathcal{B}(v_{\mathup{glo}},w)=0$ holds for any $v_{\mathup{glo}}$, then $w\in \tilde{V}$, 
where $\widetilde{V}=\CurlyBrackets*{v\in V\,|\, \pi(v)=0}$.
\end{lemma}
\begin{proof} For any $v_{\mathup{glo}}\in V_{\mathup{glo}}$, by using the fact in \cref{eq:global operator}, it is easy to obtain the first argument.  For the second argument, $\forall w\in V$, we have $\mathcal{B}(v_{\mathup{glo}},w)=0$. By using the contradiction and mimicking the proof steps in \cite{ye2024}, we can obtain that $\pi w=0$, then $w\in\tilde{V}.$ 
\end{proof} 
\section{Analysis}\label{sec:anal}
\begin{lemma}\label{lm1}
If the mesh size $H$, the wavenumber $k$ satisfies the resolution condition such that
\begin{equation}\label{eq:res}
C_1 \mu_{\mathup{msh}}^{-1}k^2H^2\leq C_*\Lambda,   
\end{equation}
where $C_1$ is the constant related to $\tilde{s}$ norm under $\texttt{T}$-coercivity in \cref{eq:bound s} and $0<C_*\ll 1$. Under the resolution condition, $\mathcal{T}_H$ is the operator on $\tilde{V}$, the continuous sesquilinear form of $\mathcal{B}$ defined in \cref{bform} is $\texttt{T}$-coercive on $ \tilde{V}$ such that  
\begin{equation}\label{lem2.2}
{\varepsilon_0} \Norm{v}_{\tilde{a}}^2\leq\mathcal{B}(v,\mathcal{T}_Hv),\quad \forall v\in \tilde{V}.
\end{equation}
where $\varepsilon_0$ is a positive constant depending on the high-contrast ratio $\Upsilon$. 
\end{lemma}
\begin{proof}
{The technique we employ consists of introducing a modified $\mathcal{T}$-operator, denoted by $\mathcal{T}_H$, defined on $\widetilde{V}$ such that for any $v \in \widetilde{V}$, one has $\mathcal{T}_H v \in \widetilde{V}$. Moreover, $\mathcal{Q}_H \colon L^2 \to V$ is the interpolation operator defined in \cref{lem:interpolation}.
 then we have $\mathcal{T}_H$ defined as follows:
\begin{align*}
  \mathcal{T}_Hv & = \mathcal{T} v -  \mathcal{Q}_H \pi \mathcal{T} v=\begin{cases}
v_1-\mathcal{Q}_H\pi v_1, & \text{ in } \Omega^+, \\
-v_2+2\mathcal{R}v_1 + \mathcal{Q}_H\pi v_2 - 2\mathcal{Q}_H \pi \mathcal{R}v_1, & \text{ in } \Omega^-.
\end{cases} \\
& = \begin{cases}
v_1,    & \text{ in } \Omega^+, \\
-v_2+2(\mathcal{R}- \mathcal{Q}_H\pi\mathcal{R})v_1=-v_2+2\mathcal{R}'_Hv_1, & \text{ in } \Omega^-.
\end{cases}
\end{align*}}
Here, the terms $\mathcal{Q}_H \pi v_1$ and $\mathcal{Q}_H \pi v_2$ vanish because $\pi v=0$ on each coarse element, and meanwhile $\mathcal{Q}_H\pi v|_{K_i}$ depends only on $\pi v|_{K_i}$. Recalling the boundness of $\mathcal{Q}_H$, we can see that $\mathcal{T}_H$ is a bounded operator that maps from $\tilde{V}$ to $\tilde{V}$. Then, let $$\mathcal{B}(v,\mathcal{T}_Hv)=\underbrace{a(v,\mathcal{T}_Hv)}_{I_1}-\underbrace{k^2(cv,\mathcal{T}_Hv)}_{I_2}.$$
For $I_1$, we already know that from \cite{ye2024},
$$I_1=a(v,\mathcal{T}_Hv)\geq(1-\sqrt{\Upsilon^{-1}}\Norm{\mathcal{R}'_H}_1)\Norm{v}_{\tilde{a}}^2,$$
{where $\norm{\mathcal{R}_H'}_1$ should be understood as a positive constant such that for all $v\in \Omega^+$,
\begin{equation*}
\RoundBrackets*{\int_{\Omega^-} \abs{\nabla \mathcal{R}_H'v}^2 \di x}^{1/2} \leq \norm{\mathcal{R}_H'}_1 \RoundBrackets*{\int_{\Omega^+} \abs{\nabla v}^2 \di x}^{1/2}.
\end{equation*}}
For $I_2$, by using of the resolution condition, we have
\begin{align*}
|I_2|&\leq k^2\mu_{\mathup{msh}}^{-1}H^2\int_{\Omega}|c|\mu_{\text{msh}}H^{-2}|v||\mathcal{T}_Hv|\dx x\leq C_1k^2\mu_{\mathup{msh}}^{-1}H^2\|v\|_{\tilde{s}}^2\\
&\leq C_1k^2\mu_{\mathup{msh}}^{-1}H^2(\frac{1}{\Lambda}\|v\|_{\tilde{a}}^2+\|\pi v\|_{\tilde{s}}^2)= C_1k^2\mu_{\mathup{msh}}^{-1}H^2\Lambda^{-1}\|v\|_{\tilde{a}}^2.
\end{align*}
We finally obtain
$$
\begin{aligned}
\mathcal{B}(v,\mathcal{T}_Hv)
\geq&\RoundBrackets*{1-\sqrt{\Upsilon^{-1}}\Norm{\mathcal{R}'_H}_1-C_1k^2\mu_{\mathup{msh}}^{-1}H^2\Lambda^{-1}}\Norm{v}_{\tilde{a}}^2\\
\geq&\RoundBrackets*{1-\sqrt{\Upsilon^{-1}}\Norm{\mathcal{R}'_H}_1-C_*}\Norm{v}_{\tilde{a}}^2\\
\geq& \varepsilon_0\norm{v}_{\tilde{a}}^2.
\end{aligned}
$$
where $\varepsilon_0$ is a positive  constant depending on $\Upsilon$ and $0<C_*\ll 1$. When we choose $\Upsilon$ to be a large value, the resolution condition \cref{eq:res} can guarantee $\varepsilon_0$ to be always positive.
\end{proof}
Based on \cref{lm1}, we now establish the well-posedness of the global problem \cref{eq:global solution}, and we need to examine the inf-sup stability of $V_{\mathup{glo}}$, i.e., proving the existence of a uniform lower bound of
\begin{equation}\label{eq:glo}
\inf_{u_{\mathup{glo}}\in V_{\mathup{glo}}} \sup_{w_{\mathup{glo}}\in V_{\mathup{glo}}}\frac{\mathcal{B}(u_{\mathup{glo}}, w_{\mathup{glo}})}{\norm{u_{\mathup{glo}}}_{\tilde{a}} \norm{w_{\mathup{glo}}}_{\tilde{a}}}.
\end{equation}
The Fortin trick \cite{Boffi2013}  suggests that it suffices to check that
\begin{equation}\label{glo}
\inf_{u\in\tilde{V}} \sup_{ w\in \tilde{V}} \frac{\mathcal{B}(u, w)}{\norm{u}_{\tilde{a}} \norm{w}_{\tilde{a}}}\geq \beta_{\text{glo}}
\end{equation}
holds. {In \cref{glo}, due to the sign-changing nature of the operator $\mathcal{B}$, the standard coercivity property no longer holds. After incorporating the $\texttt{T}$-coercivity framework via the modified transformation operator $\mathcal{T}_H$ defined in \cref{lm1}, which satisfies $\mathcal{T}_H u \in \tilde{V}$ for any $u \in \tilde{V}$, it becomes necessary to examine the coercivity of the resulting bilinear form 
$ \mathcal{B}(u, \mathcal{T}_H u)$ for  $\forall u \in \tilde{V}.$
From \cref{lm1}, for any $u \in \widetilde{V}$, it yields
\[
\mathcal{B}(u, \mathcal{T}_H u)
\ge \beta_{\text{glo}} \|u\|_{\tilde{a}}^2.
\]
which implies the desired inf-sup condition \cref{glo} with $\beta_{\text{glo}} = \varepsilon_0$.}
{\begin{theorem}\label{glothm}
Let $u$ be the solution to the original problem \cref{vf} and $u_{\mathup{glo}}$ be the solution to the global problem \eqref{eq:global solution}. If $f\in L^2(\Omega)$, then 
$$\Norm{u-u_{\mathup{glo}}}_{\tilde{a}}\leq \frac{C_1\Norm{f}_{s^{-1}}}{\varepsilon_0\sqrt{\Lambda}},$$
where the $s^{-1}$-norm is defined as
\[
  \norm{f}_{\tilde{s}^{-1}}\coloneqq \RoundBrackets*{\int_\Omega \abs{\mu}^{-1}\abs{f}^2\di x}^{1/2}
  \approx \mathcal{O}(H).
\]
\end{theorem}
\begin{proof} 
For the global error estimate in the energy norm, we can obtain that for $u-u_{\text{glo}}\in\tilde{V}$, 
\begin{equation*}
\begin{aligned}
\|u-u_{\text{glo}}\|_{\tilde{a}}^2
&\leq \frac{1}{\varepsilon_0}\mathcal{B}(u-u_{\text{glo}},\mathcal{T}_H(u-u_{\text{glo}}))
= \frac{1}{\varepsilon_0}\mathcal{B}(u,\mathcal{T}_H(u-u_{\text{glo}})) \\
&= {\frac{1}{\varepsilon_0}} (f,\mathcal{T}_H(u-u_{\text{glo}}))
\leq {\frac{C_1}{\varepsilon_0}} \|f\|_{\tilde{s}^{-1}} \|u-u_{\text{glo}}\|_{\tilde{s}} \\
&\leq \frac{{C_1}\|f\|_{\tilde{s}^{-1}}}{{\varepsilon_0}\sqrt{\Lambda}}
\|u-u_{\text{glo}}\|_{\tilde{a}}.
\end{aligned}
\end{equation*}
The first inequality follows from \cref{lm1}, which implies that $\mathcal{T}_H (u - u_{\text{glo}}) \in \tilde{V}$ for any $u - u_{\text{glo}} \in \tilde{V}$. Moreover, by \cref{lm3.2}, we have $\mathcal{B}\big(u_{\text{glo}}, \mathcal{T}_H (u - u_{\text{glo}})\big) = 0$ since $u_{\text{glo}} \in V_{\text{glo}}$ and $\mathcal{T}_H (u - u_{\text{glo}}) \in \tilde{V}$. The second inequality follows from \cref{eq:bound s}, where $C_1$ is a constant independent of $H$. Combining these results yields the desired estimate.
\end{proof}}

{
Although the global formulation yields optimal error estimates in \cref{glothm}, it is not suitable for practical computations due to its high computational cost. We therefore introduce a localized version of the multiscale method in the following \cref{eq:new local operator}. For the oversampling region $K_i^m$ defined in \cref{msbasis1}, 
the incompatibility between the standard definition of the $K_i^m$ and the \texttt{T}-coercivity framework 
necessitates the use of the symmetrization technique introduced in \cite{Chaumont2021}, which derives us to make \cref{ass:modified subdomains}.
This approach enables the construction of modified subdomains on which the exponential decay 
of the multiscale basis functions outside the corresponding symmetric patches can be established, 
as shown in \cref{lem1}.}

{
To distinguish these modified regions, we adopt the hat notation ``$\hat{\cdot}$''; 
for example, $\hat{V}_i^m = H_0^1(\hat{K}_i^m)$ denotes the space associated with the symmetric patch $\hat{K}_i^m$, 
rather than the original oversampling region $K_i^m$. 
Similarly, we establish the well-posedness of the modified local operator $\hat{\mathcal{G}}_i^m$ 
associated with the coarse element $K_i$:}
\begin{equation}\label{eq:new local operator}
\begin{split}
& \text{find } \hat{\mathcal{G}}_i^m \psi \in \hat{V}_i^m \text{ s.t. } \forall v \in \hat{V}_i^m,                                                                          \\
& \mathcal{B}(\hat{\mathcal{G}}_i^m \psi, w)_{\hat{K}_i^m} + s(\pi\hat{\mathcal{G}}_i^m \psi, \pi w)_{\hat{K}_i^m} = s(\pi\psi, \pi w)_{K_i}.
  \end{split}
\end{equation}   
Now, our multiscale finite element space $V_{\mathup{ms}}$ can be defined by solving {a} variational problem \cref{msbasis1}
$$
V_{\mathup{ms}}=\mathup{span}\CurlyBrackets*{\hat{\mathcal{G}}_i^m \psi\,|\,\, m \geq 1,  1\leq i\leq N}.
$$
Similarly, for the multiscale problem, we need to find the approximated solution of  \cref{vf}: Find $u_{\mathup{ms}}\in V_{\mathup{ms}}$ such that 
\begin{equation}\label{cemvar}
\mathcal{B}(u_{\mathup{ms}},w)=(f,w),\quad \forall w\in V_{\mathup{ms}}.
\end{equation}
\begin{assumption}\label{ass:modified subdomains}
  There exists a list of subdomains $\CurlyBrackets{\hat{K}_i^1,\hat{K}_i^2,\dots}$ that fulfills the requirements listed below:
  \begin{enumerate}
    \item Each $\hat{K}_i^m$ consists of coarse elements from $\mathcal{K}_H$.
    \item There exists an inclusion relation $K_i \subset \hat{K}_i^1 \subset \hat{K}_i^2 \subset \dots \subset$, such that $\Dist(\partial \hat{K}_i^{m},\partial \hat{K}_i^{m+1})\geq C_\mathup{msh}H$.
    \item For any $v \in V$ with $\Supp v \subset \overline{\hat{K}_i^m}$ or $\overline{\Omega\setminus \hat{K}_i^m}$, it holds that $\Supp \mathcal{T}v \subset \overline{\hat{K}_i^m}$ or $\overline{\Omega\setminus \hat{K}_i^m}$ accordingly.
  \end{enumerate}
\end{assumption}
 Before we provide the detailed proof of \cref{lem1}, we need to define the {cut-off} functions with respect to these oversampling domains as follows. 
We define a cut-off function $\hat\eta_i^{m-1,m}\in  C^{0,1}(\Omega)$ as 
\begin{equation*}
\hat{\eta}_i^{m-1,m}=
\begin{cases}
1,&\text{in}\,\,\hat{K}_i^{m-1},\\
0,&\text{in}\,\, \Omega\setminus \hat{K}_i^{m},
\end{cases}
\end{equation*}
with $0\leq \hat{\eta}_i^{m-1,m}\leq 1$ in $\hat{K}_i^{m}\setminus \hat{K}_i^{m-1}.$  By using controlling the parameter $\mu_{\text{msh}}$, we can derive the inequality
$$\abs{\nabla\hat\eta_i^{m-1,m}}^2\abs{\sigma}\leq \abs{\mu}$$
on $\Omega$ for $i$ and $m$.
\begin{assumption}\label{ass:overlapping}
  The number of coarse elements within $\hat{K}_i^m$ satisfies a relation:
  \[
    \#\CurlyBrackets*{K\in \mathcal{K}_H\mid K\subset \hat{K}_i^m} \leq C_\mathup{ol}m^d,
  \]
  where $C_\mathup{ol}$ depends solely on the mesh quality.
\end{assumption}
\begin{lemma}\label{lem1}
There exists a positive constant $\theta$ with $\theta<1$ that independents on $k,m$ and $H$ such that for any $m\geq 1$, $1\leq i\leq N_{\text{elem}}$ and $\psi\in L^2(\Omega)$,
$$\Norm{\mathcal{G}_i^{\infty}\psi}_{\tilde{a},\Omega\setminus \hat{K}_i^m}^2+\Norm{\pi\mathcal{G}_i^{\infty}\psi}_{\tilde{s},\Omega\setminus \hat{K}_i^m}^2\leq \theta^m(\Norm{\mathcal{G}_i^{\infty}\psi}_{\tilde{a}}^2+\Norm{\pi\mathcal{G}_i^{\infty}\psi}_{\tilde{s}}^2),$$
where $\theta=\RoundBrackets*{1+\RoundBrackets*{\frac{C}{2C_0+C_1+C_*\Lambda}}^2}^{-1}< 1$, {and $C=\min\{C_2-C_*,C_2-C_*\Lambda\}$.}
\end{lemma}
\begin{proof}
Taking $z_i=(1-\hat{\eta}_i^{m-1,m})\mathcal{G}_i^{\infty}\psi$, we can see that $z_i$ is supported in $\Omega\setminus \hat{K}_i^{ m-1}$. Substituting $\mathcal{T}z_i$ for $w$ in the variational form \cref{eq:global operator}, we have  
$$\mathcal{B}(\mathcal{G}_i^{\infty}\psi,\mathcal{T}z_i)+s(\pi\mathcal{G}_i^{\infty}\psi,\pi\mathcal{T}z_i)=0.$$
As a result, we obtain the following decomposition: 
\begin{align*}
&a(z_i,\mathcal{T}z_i) + s(\pi z_i,\pi \mathcal{T}z_i)  \\
&= -a\!\left(\hat{\eta}_i^{m-1,m}\mathcal{G}_i^{\infty}\psi,\mathcal{T}z_i\right)
   - s\!\left(\pi(\hat{\eta}_i^{m-1,m}\mathcal{G}_i^{\infty}\psi),\pi \mathcal{T}z_i\right)
   + k^2(c\mathcal{G}_i^{\infty}\psi,\mathcal{T}z_i) \\
&= -\underbrace{\int_{\Omega} 
    \sigma\, \mathcal{G}_i^\infty \psi \,
    \nabla \hat{\eta}_i^{m-1,m} \cdot \nabla \mathcal{T}z_i \, \mathrm{d}x}_{I_1}
   -\underbrace{\int_{\Omega} 
    \sigma\, \hat{\eta}_i^{m-1,m} 
    \nabla \mathcal{G}_i^\infty \psi \cdot \nabla \mathcal{T}z_i 
    \, \mathrm{d}x}_{I_2} \\
&\quad
   -\underbrace{\int_{\Omega}
    \mu\, \pi(\hat{\eta}_i^{m-1,m}\mathcal{G}_i^{\infty}\psi)
    \cdot \pi \mathcal{T}z_i \, \mathrm{d}x}_{I_3}
   +\underbrace{k^2 \int_{\Omega}
    c\, \mathcal{G}_i^{\infty}\psi \cdot \mathcal{T}z_i
    \, \mathrm{d}x}_{I_4}.
\end{align*}
For $I_1$, by the Cauchy-Schwarz inequality, \cref{eq:bound a} and \cref{pi1}, it is evident that
\begin{align*}
\abs{I_1} & \leq \RoundBrackets*{\int_\Omega \abs{\sigma} \abs{\nabla \hat{\eta}_{i}^{m-1,m}}^2 \abs{\mathcal{G}_i^\infty \psi}^2 \di x }^{1/2} \RoundBrackets*{\int_\Omega \abs{\sigma} \abs{\nabla \mathcal{T} z_i}^2  \di x }^{1/2}                                                               \\
& \leq \norm{\mathcal{G}_i^\infty\psi}_{\tilde{s}, \hat{K}_i^{m}\setminus \hat{K}_i^{m-1}} \norm{\mathcal{T}z_i}_{\tilde{a}, \Omega \setminus \hat{K}_i^{m-1}}           \\
& \leq C_0\RoundBrackets*{\norm{\pi\mathcal{G}_i^\infty\psi}_{\tilde{s}, \hat{K}_i^{m}\setminus \hat{K}_i^{m-1}}^2+\Lambda^{-1}\norm{\mathcal{G}_i^\infty \psi}_{\tilde{a}, \hat{K}_i^{m}\setminus \hat{K}_i^{m-1}}^2}^{1/2} \norm{z_i}_{\tilde{a}, \Omega \setminus \hat{K}_i^{m-1}}.
\end{align*}
 For $I_2$, we have
\begin{align*}
\abs{I_2} & \leq \RoundBrackets*{\int_{\hat{K}_i^{m}\setminus \hat{K}_i^{m-1}}\abs{\sigma} \abs{\hat{\eta}_i^{m-1,m}}^2 \abs{\nabla \mathcal{G}_i^\infty \psi}^2 \di x }^{1/2} \RoundBrackets*{\int_\Omega \abs{\sigma} \abs{\nabla \mathcal{T} z_i}^2  \di x }^{1/2} 
\\
& \leq \norm{\mathcal{G}_i^\infty\psi}_{\tilde{a}, \hat{K}_i^{m}\setminus \hat{K}_i^{m-1}}\norm{\mathcal{T}z_i}_{\tilde{a}, \Omega \setminus \hat{K}_i^{m-1}}                                       \leq C_0\norm{\mathcal{G}_i^\infty\psi}_{\tilde{a}, \hat{K}_i^{m}\setminus \hat{K}_i^{m-1}}\norm{z_i}_{\tilde{a}, \Omega \setminus \hat{K}_i^{m-1}}.
  \end{align*}
For $I_3$, we can similarly show that
\begin{align*}
\abs{I_3} & \leq \RoundBrackets*{\int_{\hat{K}_i^{m}\setminus \hat{K}_i^{m-1}} \abs{\mu} \abs{\pi\RoundBrackets*{\hat{\eta}_i^{m-1,m}\mathcal{G}_i^\infty\psi} }^2  \di x }^{1/2} \RoundBrackets*{\int_\Omega \abs{\mu} \abs{\pi \mathcal{T} z_i}^2  \di x }^{1/2}                                            \\
 & = \norm{\pi \RoundBrackets*{\hat{\eta}_i^{m-1,m}\mathcal{G}_i^\infty\psi}}_{\tilde{s}, \hat{K}_i^{m}\setminus \hat{K}_i^{m-1}} \norm{\pi\mathcal{T}z_i}_{\tilde{s}, \Omega \setminus \hat{K}_i^{m-1}}                                                 \\
 & \leq \norm{\hat{\eta}_i^{m-1,m}\mathcal{G}_i^\infty\psi}_{\tilde{s}, \hat{K}_i^{m}\setminus \hat{K}_i^{m-1}} \norm{\mathcal{T}z_i}_{\tilde{s}, \Omega \setminus \hat{K}_i^{m-1}}                                           \\
 & \leq C_1\norm{\mathcal{G}_i^\infty\psi}_{\tilde{s}, \hat{K}_i^{m}\setminus \hat{K}_i^{m-1}} \norm{z_i}_{\tilde{s}, \Omega \setminus \hat{K}_i^{m-1}}                                                  \\
& \leq C_1 \RoundBrackets*{\norm{\pi\mathcal{G}_i^\infty\psi}_{\tilde{s}, \hat{K}_i^{m}\setminus \hat{K}_i^{m-1}}^2+\Lambda^{-1} \norm{\mathcal{G}_i^\infty \psi}_{\tilde{a}, \hat{K}_i^{m}\setminus \hat{K}_i^{m-1}}^2}^{1/2} \RoundBrackets*{\norm{\pi z_i}_{\tilde{s}, \Omega\setminus \hat{K}_i^{m-1}}^2+\Lambda^{-1} \norm{z_i}_{\tilde{a}, \Omega \setminus \hat{K}_i^{m-1}}^2}^{1/2}.
\end{align*}
For $I_4$, by using of the resolution condition \cref{eq:res} and  \cref{eq:bound s} , we have
\begin{align*}
|I_4|\leq& k^2\mu_{\mathup{msh}}^{-1}H^2\int_{\Omega}|c|\mu_{\text{msh}}H^{-2}|z_i+\hat{\eta}_i^{m-1,m}\hat{\mathcal{G}}_i^{\infty}\psi|\cdot|\mathcal{T}z_i|\dx x\\
\leq&k^2\mu_{\mathup{msh}}^{-1}H^2\|z_i\|_{\tilde{s},\Omega\setminus \hat{K}_i^{m-1}}\|\mathcal{T}z_i\|_{\tilde{s},\Omega\setminus \hat{K}_i^{m-1}}+k^2\mu_{\mathup{msh}}^{-1}H^2\|\mathcal{G}_i^{\infty}\psi\|_{\tilde{s},\hat{K}_i^m\setminus \hat{K}_i^{m-1}}\|\mathcal{T}z_i\|_{\tilde{s},\Omega\setminus \hat{K}_i^{m-1}}\\
\leq&C_1k^2\mu_{\mathup{msh}}^{-1}H^2\|z_i\|_{\tilde{s},\Omega\setminus \hat{K}_i^{m-1}}^2+C_1k^2\mu_{\mathup{msh}}^{-1}H^2\|\mathcal{G}_i^{\infty}\psi\|_{\tilde{s},\hat{K}_i^m\setminus \hat{K}_i^{m-1}}\|z_i\|_{\tilde{s},\Omega\setminus \hat{K}_i^{m-1}}\\
\leq&C_*\Lambda\RoundBrackets*{\Lambda^{-1}\Norm{z_i}_{\tilde{a},\Omega\setminus \hat{K}_i^{m-1}}^2+\Norm{\pi z_i}_{\tilde{s},\Omega\setminus \hat{K}_i^{m-1}}^2}\\
&+C_*\Lambda \RoundBrackets*{\Lambda^{-1}\Norm{\hat{\mathcal{G}}_i^{\infty}\psi}_{\tilde{a},\hat{K}_i^m\setminus \hat{K}_i^{m-1}}^2+\Norm{\pi \hat{\mathcal{G}}_i^{\infty}\psi}_{\tilde{s},\hat{K}_i^m\setminus \hat{K}_i^{m-1}}^2}^{1/2}\|z_i\|_{\tilde{s},\Omega\setminus \hat{K}_i^{m-1}}.
\end{align*}
{Providing that $\frac{1}{\Lambda}$ is small enough ($\frac{1}{\Lambda}\leq 1)$ in \cref{eq:coer s}} and $0<C_*\ll 1$, we can obtain that
$$
\begin{aligned}
a(z_i, \mathcal{T}z_i) + s(\pi z_i, \pi\mathcal{T}z_i) 
&\geq \RoundBrackets*{\underbrace{(C_2-C_*)}_{p_1}\norm{z_i}_{\tilde{a}, \Omega\setminus \hat{K}_i^{m-1}}^2 + \underbrace{(C_2-C_*\Lambda)}_{p_2}\norm{\pi z_i}_{\tilde{s}, \Omega\setminus \hat{K}_i^{m-1}}^2}\\
&\geq C\RoundBrackets*{\norm{z_i}_{\tilde{a}, \Omega\setminus \hat{K}_i^{m-1}}^2 + \norm{\pi z_i}_{\tilde{s}, \Omega\setminus \hat{K}_i^{m-1}}^2},
\end{aligned}
$$
where $p_1\geq p_2>0$, for simplicity, we denote $C=\min\{C_2-C_*,C_2-C_*\Lambda\}$. Due to $\frac{1}{\Lambda}$ is small enough ($\frac{1}{\Lambda}\leq 1)$, we have 
\begin{align*}
& \quad \RoundBrackets*{\norm{\mathcal{G}_i^\infty \psi}_{\tilde{a}, \Omega\setminus K_i^{m}}^2 + \norm{\pi \mathcal{G}_i^\infty \psi}_{\tilde{s}, \Omega\setminus K_i^{m}}^2}^{1/2} \leq \RoundBrackets*{\norm{z_i}_{\tilde{a}, \Omega\setminus \hat{K}_i^{m-1}}^2 + \norm{\pi z_i}_{\tilde{s}, \Omega\setminus \hat{K}_i^{m-1}}^2}^{1/2} \\
& \leq \frac{2C_0+C_1+C_*\Lambda}{C} \RoundBrackets*{\norm{\mathcal{G}_i^\infty \psi}_{\tilde{a}, \hat{K}_i^{m}\setminus \hat{K}_i^{m-1}}^2+\norm{\pi\mathcal{G}_i^\infty\psi}_{\tilde{s}, \hat{K}_i^{m}\setminus {\hat{K}}_i^{m-1}}^2}^{1/2}.
\end{align*}
We hence derive an iteration relation
\begin{align*}
& \norm{\mathcal{G}_i^\infty \psi}_{\tilde{a}, \Omega\setminus \hat{K}_i^{m}}^2+\norm{\pi\mathcal{G}_i^\infty\psi}_{\tilde{s}, \Omega \setminus\hat{K}_i^{m}}^2 \\
& \qquad \leq  \RoundBrackets*{1+\RoundBrackets*{\frac{C}{2C_0+C_1+C_*\Lambda}}^2}^{-1} \RoundBrackets*{\norm{\mathcal{G}_i^\infty \psi}_{\tilde{a}, \Omega \setminus \hat{K}_i^{m-1}}^2+\norm{\pi\mathcal{G}_i^\infty\psi}_{\tilde{s}, \Omega\setminus \hat{K}_i^{m-1}}^2}
,
\end{align*}
  which completes the proof.
\end{proof}
\begin{lemma}\label{lem2}
It holds that for any $m\geq 1, 1\leq i\leq N$ and $\psi\in L^2(\Omega)$,
$$\norm{(\mathcal{G}_i^\infty-\hat{\mathcal{G}}_i^m) \psi}_{\tilde{a}}^2+\norm{\pi(\mathcal{G}_i^\infty-\hat{\mathcal{G}}_i^m) \psi}_{\tilde{s}}^2\leq C'\theta^{m-1}\RoundBrackets*{\norm{\mathcal{G}_i^\infty \psi}_{\tilde{a}}^2+\norm{\pi\mathcal{G}_i^\infty \psi}_{\tilde{s}}^2},$$
where  $0<\theta< 1$ comes from \cref{lem1} and $C'=\RoundBrackets*{\frac{4C_0+C_1+C_*\Lambda}{C}}^2>0$ depends on $C,C_0,C_1,C_*$ and $\Lambda$.
\end{lemma}
\begin{proof}
We now take $z_i=(\mathcal{G}_i^\infty-\hat{\mathcal{G}}_i^m)\psi$ and  introduce a decomposition 
$$z_i=\underbrace{(1-\hat{\eta}_i^{m-1,m})\mathcal{G}_i^\infty\psi}_{z_i'}+\underbrace{(\hat{\eta}_i^{m-1,m}-1)\hat{\mathcal{G}}_i^m\psi+\hat{\eta}_i^{m-1,m}z_i}_{z_i^{''}}.$$   
We can observe that $z_i^{''}\in \hat{V}_i^m$, which implies $\mathcal{T}z_i^{''}\in \hat{V}_i^m$. Therefore, it is easy to see that $\mathcal{B}(z_i,\mathcal{T}z_i^{''})+s(\pi z_i,\pi\mathcal{T}z_i^{''})=0$. Now we start to estimate $\mathcal{B}(z_i,\mathcal{T}z_i')+s(\pi z_i,\pi\mathcal{T}z_i')$, and the next step is 
\begin{align*}
a(z_i,\mathcal{T}z_i)+s(\pi z_i,\mathcal{T}z_i)&=\mathcal{B}(z_i,\mathcal{T}z_i)+s(\pi z_i,\mathcal{T}z_i)+k^2(z_i,\mathcal{T}z_i)\\
&=\mathcal{B}(z_i,\mathcal{T}z_i')+s(\pi z_i,\pi\mathcal{T}z_i')+k^2(z_i,\mathcal{T}z_i)\\
&=a(z_i,\mathcal{T}z_i')+s(\pi z_i,\pi\mathcal{T}z_i')+k^2(z_i,\mathcal{T}z_i'')\\
&\leq\underbrace{\|z_i\|_{\tilde{a}}\|\mathcal{T}z_i'\|_{\tilde{a}}}_{I_1}+\underbrace{k^2(z_i,\mathcal{T}z_i^{''})}_{I_2}+\underbrace{\|\pi z_i\|_{\tilde{s}}\|\pi\mathcal{T}z_i'\|_{\tilde{s}}}_{I_3}.
\end{align*}
For $I_1$, we recall that $\|\mathcal{T}z_i'\|_{\tilde{a}}\leq C_0\|z_i'\|_{\tilde{a}}$, and for $z_i'$ we have
\begin{align*}
\norm{z'_i}_{\tilde{a}}^2 & \leq 2\int_{\Omega \setminus {\hat{K}}_i^{m-1}}  \abs{\sigma}\abs{1-\hat{\eta}_i^{m-1,m}}^2 \abs{ \nabla {\hat{\mathcal{G}}}_i^\infty \psi}^2 \di x + 2 \int_{\hat{K}_i^{m}\setminus {\hat {K}}_i^{m-1}} \abs{\sigma}\abs{\nabla \hat{\eta}_i^{m-1,m}}^2 \abs{{\hat{\mathcal{G}}}_i^\infty\psi}^2 \di x \\
& \leq 2 \norm{{\mathcal{G}}_i^\infty \psi}^2_{\tilde{a},\Omega \setminus {\hat{K}}_i^{m-1}} + 2 \int_{\Omega\setminus {\hat{K}}_i^{m-1}} \abs{\mu} \abs{\mathcal{G}_i^\infty \psi}^2 \di x\\
& \leq 2(1+\Lambda^{-1})\norm{{\mathcal{G}}_i^\infty \psi}^2_{\tilde{a},\Omega \setminus {\hat{K}}_i^{m-1}} + 4 \norm{\pi \mathcal{G}_i^\infty \psi}_{\tilde{s}, \Omega \setminus {\hat{K}}_i^{m-1}}^2.
\end{align*}
For $I_2$, we use the resolution condition \cref{eq:res} and the definition of $z_i^{''}$ to obtain 
\begin{align*}
|I_2|&\leq k^2\mu_{\mathup{msh}}^{-1}H^2\int_{\Omega}|c|\mu_{\text{msh}}H^{-2}|z_i|\cdot|\mathcal{T}z_i''|\dx x\leq C_1k^2\mu_{\mathup{msh}}^{-1}H^2\|z_i\|_{\tilde{s}}\|z_i''\|_{\tilde{s}}\\
&\leq C_*\Lambda\RoundBrackets*{\|z_i\|_{\tilde{s}}^2+\|z_i\|_{\tilde{s}}\|z_i'\|_{\tilde{s}}}\\
&\leq C_*\Lambda\RoundBrackets*{\Lambda^{-1}\|z_i\|_{\tilde{a}}^2+\|\pi z_i\|_{\tilde{s}}^2}+C_*\Lambda\RoundBrackets*{\Lambda^{-1}\| z_i\|_{\tilde{a}}^2+\| \pi z_i\|_{\tilde{s}}^2}^{1/2}\|z_i'\|_{\tilde{s}}.
\end{align*}
For $I_3$, $\norm{\pi\mathcal{T}z'_i}_{\tilde{s}} \leq \norm{\mathcal{T}z'_i}_{\tilde{s}} \leq C_1\norm{z'_i}_{\tilde{s}} $ and
\begin{align*}
\norm{z'_i}_{\tilde{s}}^2 & =\int_{\Omega} \abs{\mu} \abs{1-\hat{\eta}_i^{m-1,m}}^2 \abs{\mathcal{G}_i^\infty \psi}^2 \di x \leq \norm{\mathcal{G}_i^\infty \psi}_{\tilde{s},\Omega\setminus \hat{K}_i^{m-1}}^2  \\
&\leq \Lambda^{-1} \norm{\mathcal{G}_i^\infty \psi}_{\tilde{a},\Omega\setminus \hat{K}_i^{m-1}}^2 + \norm{\pi\mathcal{G}_i^\infty \psi}_{\tilde{s}, \Omega\setminus \hat{K}_i^{m-1}}^2.
\end{align*}
Then, providing that $\frac{1}{\Lambda}$ and $C_*$ {are} small enough ($\frac{1}{\Lambda},C_*\leq 1)$, we have the following estimate
\begin{align*}
C\RoundBrackets*{\norm{z_i}_{\tilde{a}}^2 +\norm{\pi z_i}_{\tilde{s}}^2}  \leq a(z_i, \mathcal{T}z_i) + s(\pi z_i, \pi\mathcal{T}z_i).
\end{align*}
We { can finally  obtain that}
  \[
    \norm{z_i}_{\tilde{a}}^2 + \norm{\pi z_i}_{\tilde{s}}^2 \leq \RoundBrackets*{\frac{4C_0+C_1+C_*\Lambda}{C}}^2 \RoundBrackets*{\norm{\mathcal{G}_i^\infty \psi}^2_{\tilde{a},\Omega\setminus \hat{K}_i^{m-1}}+\norm{\pi \mathcal{G}_i^\infty \psi}_{\tilde{s},\Omega\setminus \hat{K}_i^{m-1}}^2},
    \]
{and} we hence derive the desired result by utilizing \cref{lem1}.
\end{proof}
\begin{theorem}\label{prop:exponential decay}
Under \cref{ass:overlapping}, there exists constants $\Upsilon,\Lambda^*$ such that for any $\Upsilon\geq \Upsilon'$ and $\Lambda\geq \Lambda^*$, it holds that for any $m\geq 1, 1\leq i\leq N$ and any $\psi \in L^2(\Omega)$,
\[
\norm{(\mathcal{G}^\infty-\hat{\mathcal{G}}^m) \psi}_{\tilde{a}}^2 + \norm{\pi (\mathcal{G}^\infty-\hat{\mathcal{G}}^m)\psi}_{\tilde{s}}^2 \leq c_0(\Lambda)C_\mathup{ol} (m+1)^d \theta^{m-1}\norm{\pi\psi}_{\tilde{s}}^2,
  \]
where \(0 < \theta < 1\) is given in \cref{lem1}, and \(c_0(\Lambda) > 0\) is independent of \(H\).
\end{theorem}
\begin{proof}
 We denote $z_i=(\mathcal{G}_i^{\infty}-\hat{\mathcal{G}}_i^m)\psi$ and $z=\sum_{i=1}^{N}z_i$. Then we take the decomposition of $z$ as
 $$z=\underbrace{(1-\hat{\eta}_i^{m,m+1})z}_{z'}+\underbrace{\hat{\eta}_i^{m.m+1}z}_{z''}$$
 where $i$ is arbitrary chosen from $1,\dots,N$. We can observe that $\mathup{supp}(\mathcal{T}z')\cap K_i=\emptyset $, then we have $\mathcal{B}(z_i,\mathcal{T}z')+s(\pi z_i,\pi\mathcal{T}z')=0$. 
 The next step
 \begin{align*}
a(z_i,\mathcal{T}z)+s(\pi z_i,\pi\mathcal{T}z)&=\mathcal{B}(z_i,\mathcal{T}z'')+s(\pi z_i,\pi\mathcal{T}z'')+k^2(z_i,\mathcal{T}z)\\
&=a(z_i,\mathcal{T}z'')+s(\pi z_i,\pi\mathcal{T}z'')+k^2(z_i,\mathcal{T}z')\\
&\leq\underbrace{\|z_i\|_{\tilde{a}}\|\mathcal{T}z''\|_{\tilde{a}}}_{I_1}+\underbrace{k^2(z_i,\mathcal{T}z')}_{I_2}+\underbrace{\|\pi z_i\|_{\tilde{s}}\|\pi\mathcal{T}z''\|_{\tilde{s}}}_{I_3}.
\end{align*}
Similarly, for $\norm{\mathcal{T}z''}_{\tilde{a}} $ in $I_1$, by using \cref{eq:coer a}
\[
\begin{aligned}
\norm{\mathcal{T}z''}_{\tilde{a}} \leq C_0 \norm{\hat{\eta}_i^{m,m+1}z}_{\tilde{a}}&\leq C_0\RoundBrackets*{\norm{z}_{\tilde{a},\hat{K}_i^{m+1}}+\norm{z}_{\tilde{s},\hat{K}_i^{m+1}}}\\
&\leq C_0\RoundBrackets*{(1+\frac{1}{\sqrt{\Lambda}})\norm{z}_{\tilde{a},\hat{K}_i^{m+1}} + \norm{\pi z}_{\tilde{s},\hat{K}_i^{m+1}}}.
\end{aligned}
\]
For $\|\pi\mathcal{T}z''\|_{\tilde{s}}$ in $I_3$, by using \cref{eq:bound s}
\[
\begin{aligned}
\norm{\pi\mathcal{T}z''}_{\tilde{s}}\leq \norm{\mathcal{T}z''}_{\tilde{s}} \leq C_1 \norm{\hat{\eta}_i^{m,m+1}z}_{\tilde{s}} 
&\leq C_1\norm{z}_{\tilde{s},\hat{K}_i^{m+1}} \\&\leq C_1\RoundBrackets*{\frac{1}{\sqrt{\Lambda}}\norm{z}_{\tilde{a},\hat{K}_i^{m+1}}+\norm{\pi z}_{\tilde{s},\hat{K}_i^{m+1}}}.
\end{aligned}
\]
For $I_2$, we use the resolution condition \cref{eq:res} and inequality \cref{eq:bound s}, 
\begin{align*}
|I_2|&\leq k^2\mu_{\mathup{msh}}^{-1}H^2\int_{\Omega}|c|\mu_{\text{msh}}H^{-2}|z_i|\cdot|\mathcal{T}z'|\dx x\leq C_1k^2\mu_{\mathup{msh}}^{-1}H^2\|z_i\|_{\tilde{s}}\|z'\|_{\tilde{s}}\\
&\leq C_*\Lambda\RoundBrackets*{
\norm{z}_{\tilde{s}}^2+\norm{z_i}_{\tilde{s}}\norm{z}_{\tilde{s},\hat{K}_i^{m+1}}}\\
&\leq C_*\Lambda\RoundBrackets*{
\Lambda^{-1}\norm{z}_{\tilde{a}}^2+\norm{\pi z}_{\tilde{s}}^2+\norm{z_i}_{\tilde{s}}\norm{z}_{\tilde{s},\hat{K}_i^{m+1}}}.
\end{align*}
To simplify the expression, we assume that $\frac{1}{\Lambda}\leq 1$, and combine the estimates of $I_1, I_2$ and $I_3$ to derive that
\begin{align*}
a(z_i, \mathcal{T}z)+s(\pi z_i,\pi\mathcal{T}z)
\leq&\underbrace{(2C_0+C_1+C_*\Lambda)}_{c_0(\Lambda)}\RoundBrackets*{\norm{z}_{\tilde{a}, \hat{K}_i^{m+1}}^2+\norm{\pi z}_{\tilde{s},\hat{K}_i^{m+1}}^2}^{1/2} \RoundBrackets*{\norm{z_i}_{\tilde{a}}^2+\norm{\pi z_i}_{\tilde{s}}^2}^{1/2}\\
&+C_*\norm{z}_{\tilde{a}}^2+C_*\Lambda\norm{\pi z}_{\tilde{s}}^2,
\end{align*}
where $c_0(\Lambda)$ depends on $\Lambda$, $C_0$, $C_1$, and $C_2$.
Based on the oversampling condition \cref{ass:overlapping}, it holds that
\[
\sum_{i=1}^{N} \RoundBrackets*{\norm{z}_{\tilde{a},\hat{K}_i^{m+1}}^2+\norm{\pi z}_{\tilde{s},\hat{K}_i^{m+1}}^2} \leq C_\mathup{ol}(m+1)^d \RoundBrackets*{ \norm{z}_{\tilde{a}}^2+\norm{\pi z}_{\tilde{s}}^2}.
  \]
Therefore, applying the Cauchy--Schwarz inequality,
\begin{align*}
& \quad C''\RoundBrackets*{\norm{z}_{\tilde{a}}^2+\norm{\pi z}_{\tilde{s}}^2} \leq \sum_{i=1}^{N} \RoundBrackets*{a(z_i, \mathcal{T}z)+s(\pi z_i,\pi\mathcal{T}z)}   \\
& \leq c_0(\Lambda) \RoundBrackets*{\sum_{i=1}^{N} \RoundBrackets*{\norm{z}_{\tilde{a},\hat{K}_i^{m+1}}^2+\norm{\pi z}_{\tilde{s},\hat{K}_i^{m+1}}^2}}^{1/2} \RoundBrackets*{\sum_{i=1}^{N}\RoundBrackets*{\norm{z_i}_{\tilde{a}}^2+\norm{\pi z_i}_{\tilde{s}}^2}}^{1/2} \\
& \leq c_0(\Lambda)\sqrt{C_\mathup{ol}} (m+1)^{d/2}\RoundBrackets*{ \norm{z}_{\tilde{a}}^2+\norm{\pi z}_{\tilde{s}}^2}^{1/2} \RoundBrackets*{\sum_{i=1}^{N}\RoundBrackets*{\norm{z_i}_{\tilde{a}}^2+\norm{\pi z_i}_{\tilde{s}}^2}}^{1/2},
\end{align*}
where $C''=\max\{C_2-C_*N,C_2-C_*\Lambda N\}>0$.  It has been shown in \cref{lem2} that, for any $i$,
\[
\norm{z_i}_{\tilde{a}}^2+\norm{\pi z_i}_{\tilde{s}}^2 \leq  C'\theta^{m-1}\RoundBrackets*{\norm{\mathcal{G}_i^\infty \psi}_{\tilde{a}}^2 + \norm{\pi\mathcal{G}_i^\infty \psi}_{\tilde{s}}^2},
\]
and we turn to provide a bound on $\sum_{i=1}^{N} \norm{\mathcal{G}_i^\infty \psi}_{\tilde{a}}^2 + \norm{\pi\mathcal{G}_i^\infty \psi}_{\tilde{s}}^2$ by letting $w=\mathcal{G}_i^{\infty}\psi$ in \cref{eq:global operator} to obtain 
\begin{align*}
& \quad  C_2 \RoundBrackets*{\norm{\mathcal{G}_i^\infty \psi}_{\tilde{a}}^2 + \norm{\pi\mathcal{G}_i^\infty \psi}_{\tilde{s}}^2} \leq a(\mathcal{G}_i^\infty\psi, \mathcal{T}\mathcal{G}_i^\infty\psi)+s(\pi\mathcal{G}_i^\infty\psi, \pi\mathcal{T} \mathcal{G}_i^\infty\psi)   \\
& = s(\pi\psi, \pi\mathcal{T}\mathcal{G}_i^\infty \psi)_{K_i} +k^2(\mathcal{G}_i^\infty\psi, \mathcal{T}\mathcal{G}_i^\infty\psi)\\
&\leq \norm{\pi\psi}_{\tilde{s},K_i}\norm{\pi\mathcal{T}\mathcal{G}_i^\infty\psi}_{\tilde{s},K_i} +\underbrace{k^2(\mathcal{G}_i^\infty\psi, \mathcal{T}\mathcal{G}_i^\infty\psi)}_{I_1} \\
& \leq \norm{\pi\psi}_{\tilde{s},K_i}\norm{\mathcal{T}\mathcal{G}_i^\infty\psi}_{\tilde{s}}+C_*\Lambda\RoundBrackets*{\Lambda^{-1}\|\mathcal{G}_i^\infty\psi\|_{\tilde{a}}^2+\|\pi\mathcal{G}_i^\infty\psi\|_{\tilde{s}}^2}\\
&\leq C_1\norm{\pi\psi}_{\tilde{s},K_i}\norm{\mathcal{G}_i^\infty\psi}_{\tilde{s}}+C_*\Lambda\RoundBrackets*{\Lambda^{-1}\|\mathcal{G}_i^\infty\psi\|_{\tilde{a}}^2+\|\pi\mathcal{G}_i^\infty\psi\|_{\tilde{s}}^2}\\
& \leq C_1\norm{\pi\psi}_{\tilde{s},K_i}\RoundBrackets*{\Lambda^{-1}\norm{\mathcal{G}_i^\infty\psi}_{\tilde{a}}^2+\norm{\pi\mathcal{G}_i^\infty\psi}_{\tilde{s}}^2}^{1/2}+C_*\Lambda\RoundBrackets*{\Lambda^{-1}\|\mathcal{G}_i^\infty\psi\|_{\tilde{a}}^2+\|\pi\mathcal{G}_i^\infty\psi\|_{\tilde{s}}^2},
\end{align*}
where for the boundedness of $I_1$, we use the resolution condition \cref{eq:res} and \cref{eq:bound s}
\begin{align*}
|I_1|&\leq k^2\mu_{\mathup{msh}}^{-1}H^2\int_{\Omega}|c|\mu_{\text{msh}}H^{-2}|\mathcal{G}_i^\infty\psi|\cdot| \mathcal{T}\mathcal{G}_i^\infty\psi|\dx x\\
&\leq C_1k^2\mu_{\mathup{msh}}^{-1}H^2\|\mathcal{G}_i^\infty\psi\|_{\tilde{s}}^2\leq C_*\Lambda\RoundBrackets*{\Lambda^{-1}\|\mathcal{G}_i^\infty\psi\|_{\tilde{a}}^2+\|\pi\mathcal{G}_i^\infty\psi\|_{\tilde{s}}^2}.
\end{align*}
Finally, the above  gives that
\begin{align*}
C(\norm{\mathcal{G}_i^\infty \psi}_{\tilde{a}}^2+\norm{\pi\mathcal{G}_i^\infty \psi}_{\tilde{s}}^2)&\leq(C_2-C_*)\norm{\mathcal{G}_i^\infty \psi}_{\tilde{a}}^2 +(C_2-C_*\Lambda)\norm{\pi\mathcal{G}_i^\infty \psi}_{\tilde{s}}^2 \leq c_1 \norm{\pi\psi}_{\tilde{s},K_i}^2\\
\norm{\mathcal{G}_i^\infty \psi}_{\tilde{a}}^2+\norm{\pi\mathcal{G}_i^\infty \psi}_{\tilde{s}}^2&\leq c_1 \norm{\pi\psi}_{\tilde{s},K_i}^2,
\end{align*}
where $c_1=(C_1/C)^2>0$. We finally completed the proof by collecting all the estimates obtained above.
\end{proof}
\begin{lemma}[ref.\ \cite{Chung2018}]\label{lem:interpolation}
  There exists a bounded map $\mathcal{Q}_H\colon L^2(\Omega) \rightarrow V$ and a positive constant $C_\mathup{inv}$ such that for all $v \in L^2(\Omega)$, it holds that $\pi\mathcal{Q}_H v=\pi v$ and $\norm{\mathcal{Q}_H v}_{\tilde{a}} \leq C_\mathup{inv} \norm{\pi v}_{\tilde{s}}$.
  Moreover, for each coarse element $K_i$, $\mathcal{Q}_Hv|_{K_i}$ depends only on the data of $v$ in $K_i$ and vanishes on $\partial K_i$.
\end{lemma}
The remainder of this section is dedicated to proving an error estimate for the multiscale solution, as stated in \cref{thm:main local} below. Before we prove the \cref{thm:main local}, we need to prove Cea's lemma of the multiscale problem \cref{cemvar}.
\begin{lemma}[Cea's Lemma]\label{cea}
If $u$ is the real solution of \cref{vf} and $u_{\mathup{ms}}$ is the multiscale solution of \cref{cemvar}, there exists a $v_{\mathup{ms}}\in V_{\mathup{ms}}$, such that
$$\|u-u_{\mathup{ms}}\|_{\tilde{a}}\leq \frac{\alpha}{\beta_{\mathup{cem}}}\|u-v_{\mathup{ms}}\|_{\tilde{a}},$$
where  {$\alpha=C_1 + C_*+C_1k^2\mu_{\mathup{msh}}^{-1}C_{\mathup{po}}^2$} and $\beta_{\mathup{cem}}>0$ is a constant from the inf-sup condition of the multiscale problem.
\end{lemma}
\begin{proof}
We only need to check the ellipticity and continuity of the bilinear form $\mathcal{B}$ over $V$, where $V=V_1\bigoplus V_{\mathup{ms}}$. To verify the continuity of $\mathcal{B}$, we consider for all $u, v \in V$ and write
\begin{align*}
    \mathcal{B}(u,\mathcal{T}v)=\underbrace{a(u,\mathcal{T}v)}_{I_1}-\underbrace{k^2(cu,\mathcal{T}v)}_{I_2}, 
\end{align*}
and for the second term $I_2$, we use the resolution condition \cref{eq:res} and \cref{eq:bound s}
{\begin{align*}
|I_2|&\leq k^2\mu_{\mathup{msh}}^{-1}H^2\int_{\Omega}|c|\mu_{\text{msh}}H^{-2}|u|\cdot|\mathcal{T}v|\dx x\\
&\leq C_0k^2\mu_{\mathup{msh}}^{-1}H^2\|u\|_{\tilde{s}}\|v\|_{\tilde{s}}\leq C_1k^2\mu_{\mathup{msh}}^{-1}H^2\RoundBrackets*{\Lambda^{-1}\|u\|_{\tilde{a}}^2+\|\pi u\|_{\tilde{s}}^2}^{1/2}\RoundBrackets*{\Lambda^{-1}\|v\|_{\tilde{a}}^2+\|\pi v\|_{\tilde{s}}^2}^{1/2}.
\end{align*}}
Meanwhile, utilizing Poincar$\acute{\mathup{e}}$ inequality, we can see that
$$\|\pi u\|_{\tilde{s}}\leq\|u\|_{\tilde{s}}\leq C_{\mathup{po}}H^{-1}\|u\|_{\tilde{a}},\quad \|\pi v\|_{\tilde{s}}\leq\|v\|_{\tilde{s}}\leq C_{\mathup{po}}H^{-1}\|v\|_{\tilde{a}}.$$
{Then we can obtain
\begin{align*}
|\mathcal{B}(u,\mathcal{T}v)|&\leq C_0\|u\|_{\tilde{a}}\|v\|_{\tilde{a}}+C_1k^2\mu_{\mathup{msh}}^{-1}H^2(\Lambda^{-1}+C_{\mathup{po}}^2H^{-2})\|u\|_{\tilde{a}}\|v\|_{\tilde{a}}\\
&\leq (C_0 + C_*+C_1k^2\mu_{\mathup{msh}}^{-1}C_{\mathup{po}}^2)\|u\|_{\tilde{a}}\|v\|_{\tilde{a}},
\end{align*}
where $\alpha = C_0 + C_* + C_1 k^2 \mu_{\mathup{msh}}^{-1} C_{\mathup{po}}^2$ is a constant independent of $H$.} Then, for the ellipticity of $\mathcal{B}$, we need to prove the inf-sup condition of the multiscale problem \cref{cemvar}, demonstrating that 
\[
  \inf_{u_{\mathup{ms}}\in V_{\mathup{ms}}} \sup_{v_{\mathup{ms}}\in V_{\mathup{ms}}} \frac{\mathcal{B}(u_{\mathup{ms}}, v_{\mathup{ms}})}{\norm{u_{\mathup{ms}}}_{\tilde{a}} \norm{v_{\mathup{ms}}}_{\tilde{a}}}\geq \beta_{\mathup{cem}},
\]
where $\beta_{\mathup{cem}}$ is a positive constant independent of $H$. First, for any $u_{\mathup{ms}}\in V_{\mathup{ms}}$, we can find $\psi$ such that $u\in V_{\mathup{ms}}=\hat{\mathcal{G}}^m\psi$. Then we choose $u_{\mathup{glo}}=\mathcal{G}^{\infty}\psi.$ Based on the inf-sup condition of the global version \cref{eq:glo}, we can find $v_{\mathup{glo}}=\mathcal{G}^{\infty}\phi\in V_{\mathup{glo}}$, such that $\mathcal{B}(u_{\mathup{glo}}, v_{\mathup{glo}})\geq\beta_{\mathup{glo}}\|u_{\mathup{glo}}\|_{\tilde{a}}\|v_{\mathup{glo}}\|_{\tilde{a}}$. If we define $v_{\mathup{ms}}=\hat{\mathcal{G}}^m\phi$, then we can obtain
\begin{align*}
  \mathcal{B}(u_{\mathup{ms}},v_{\mathup{ms}})&=\mathcal{B}(u_{\mathup{glo}},v_{\mathup{glo}})+\mathcal{B}(u_{\mathup{glo}},v_{\mathup{ms}}-v_{\mathup{glo}})+\mathcal{B}(u_{\mathup{ms}}-u_{\mathup{glo}},v_{\mathup{ms}})\\
 &\geq \beta_{\mathup{glo}}\|u_{\mathup{glo}}\|_{\tilde{a}}\|v_{\mathup{glo}}\|_{\tilde{a}}-\underbrace{\alpha\|u_{\mathup{glo}}\|_{\tilde{a}}\|v_{\mathup{ms}}-v_{\mathup{glo}}\|_{\tilde{a}}}_{I_1}-\underbrace{\alpha \|u_{\mathup{ms}}-u_{\mathup{glo}}\|_{\tilde{a}}\|v_{\mathup{ms}}\|_{\tilde{a}}}_{I_2}.
\end{align*}
For $I_1$, we need to give the estimate of $\|v_{\mathup{ms}}-v_{\mathup{glo}}\|_{\tilde{a}}=\|(\hat{\mathcal{G}}^m-\mathcal{G}^{\infty})\phi\|_{\tilde{a}}.$   {From the exponential decay properties in \cref{prop:exponential decay}, it shows that $\|(\mathcal{G}^m-\mathcal{G}^{\infty})\phi\|_{\tilde{a}}\leq C_3\theta^{m-1}\|\pi\phi\|_{\tilde{s}}$, where $C_3$ is a constant independent of $H$.} Then, the remaining task is to obtain $\|\pi\phi\|_{\tilde{s}}$ can be bounded by $\|\mathcal{G}^{\infty}\phi\|_{\tilde{a}}$.
It is not harmful to assume that $\phi \in V$ with $\norm{\phi}_{\tilde{a}} \leq C_\mathup{inv} \norm{\pi\phi}_{\tilde{s}}$ from \cref{lem:interpolation}.
Taking $\mathcal{T}\phi$ as a test function in \cref{eq:global operator}, we can obtain that
\[
  \mathcal{B}(\mathcal{G}^\infty\phi, \mathcal{T} \phi)+ s(\pi\mathcal{G}^\infty\phi, \pi\mathcal{T}\phi)=s(\pi\phi, \pi\mathcal{T}\phi).
\]
The estimate \cref{eq:coer s}  says that
\[
  C_3\norm{\pi\phi}_{\tilde{s}}^2-C_4\norm{\phi}_{\tilde{a}}^2 \leq s(\pi\phi, \pi\mathcal{T} \phi)
\]
where $C_3 \rightarrow 1$ and $C_4 \rightarrow 0$ as $\Upsilon \rightarrow \infty$.
Therefore, if $\Upsilon$ is large enough and $\Lambda^{-1}\leq 1$, we can show that
\begin{align*}
C_5\norm{\pi \phi}_{\tilde{s}}^2 & \leq  \RoundBrackets*{C_3-C_4C_\mathup{inv}^2}\norm{\pi \phi}_{\tilde{s}}^2 \leq C_3\norm{\pi \phi}_{\tilde{s}}^2-C_4\norm{\phi}_{\tilde{a}}^2 \\
& \leq s(\pi\phi, \pi\mathcal{T} \phi) \leq \alpha\norm{\mathcal{G}^\infty\phi}_{\tilde{a}} \norm{\phi}_{\tilde{a}}+\norm{\pi\mathcal{G}^\infty\phi}_{\tilde{s}} \norm{\pi\mathcal{T}\phi}_{\tilde{s}}    \\
& \leq \alpha \norm{\mathcal{G}^\infty\phi}_{\tilde{a}} \norm{\phi}_{\tilde{a}}+C_1\norm{\pi\mathcal{G}^\infty\phi}_{\tilde{s}} \norm{\phi}_{\tilde{s}}\\
& \leq \alpha \norm{\mathcal{G}^\infty\phi}_{\tilde{a}} \norm{\phi}_{\tilde{a}}+C_1\norm{\pi\mathcal{G}^\infty\phi}_{\tilde{s}} \RoundBrackets*{\Lambda^{-1}\norm{\phi}_{\tilde{a}}^2+\norm{\pi\phi}_{\tilde{s}}^2}^{1/2}\\\
& \leq \sqrt{C_1^2+\alpha^2}\RoundBrackets*{\norm{\mathcal{G}^\infty\phi}_{\tilde{a}}^2+\norm{\pi\mathcal{G}^\infty\phi}_{\tilde{s}}^2}^{1/2}\RoundBrackets*{\norm{\phi}_{\tilde{a}}^2+\norm{\pi\phi}_{\tilde{s}}^2}^{1/2}\\
& \leq \sqrt{C_1^2+\alpha^2}\sqrt{1+C_{\mathup{inv}}^2}\RoundBrackets*{\norm{\mathcal{G}^\infty\phi}_{\tilde{a}}^2+\norm{\pi\mathcal{G}^\infty\phi}_{\tilde{s}}^2}^{1/2}\norm{\pi\phi}_{\tilde{s}}.
\end{align*}
Finally, after using Poincar$\acute{\mathup{e}}$ inequality, we can see that
$$\|\pi\mathcal{G}^{\infty}\phi\|_{\tilde{s}}\leq\|\mathcal{G}^{\infty}\phi\|_{\tilde{s}}\leq C_{\mathup{po}}H^{-1}\|\mathcal{G}^{\infty}\phi\|_{\tilde{a}}.$$
{Then, $\|v_{\mathup{ms}} - v_{\mathup{glo}}\|_{\tilde{a}} \le C_*' C_3 \theta^{m-1} H^{-1} \|\mathcal{G}^{\infty}\phi\|_{\tilde{a}} = \eta \|v_{\mathup{glo}}\|_{\tilde{a}},$ where $C_*' = \frac{2C_{\mathup{po}}\sqrt{C_1^2+\alpha^2}\sqrt{1+C_{\mathup{inv}}^2}\,}{C_5}$, and $\eta > 0$ is sufficiently small due to the exponential decay factor $\theta^{m-1}$.}
Following the techniques in $I_1$ for $I_2$, we can obtain$\|u_{\mathup{ms}}-u_{\mathup{glo}}\|_{\tilde{a}}\leq C_*' C_3 \theta^{m-1}H^{-1}\|\mathcal{G}^{\infty}\psi\|_{\tilde{a}}=\eta\|u_{\mathup{glo}}\|_{\tilde{a}}.$ Combining with all together, we can obtain the inf-sup condition 
\begin{align*}
\mathcal{B}(u_{\mathup{ms}},v_{\mathup{ms}})&=\mathcal{B}(u_{\mathup{glo}},v_{\mathup{glo}})+\mathcal{B}(u_{\mathup{glo}},v_{\mathup{ms}}-v_{\mathup{glo}})+\mathcal{B}(u_{\mathup{ms}}-u_{\mathup{glo}},v_{\mathup{ms}})\\
&\geq \beta_{\mathup{glo}}\|u_{\mathup{glo}}\|_{\tilde{a}}\|v_{\mathup{glo}}\|_{\tilde{a}}-\underbrace{\alpha\|u_{\mathup{glo}}\|_{\tilde{a}}\|v_{\mathup{ms}}-v_{\mathup{glo}}\|_{\tilde{a}}}_{I_1}-\underbrace{\alpha \|u_{\mathup{ms}}-u_{\mathup{glo}}\|_{\tilde{a}}\|v_{\mathup{ms}}\|_{\tilde{a}}}_{I_2}\\
&\geq \beta_{\mathup{glo}}\|u_{\mathup{glo}}\|_{\tilde{a}}\|v_{\mathup{glo}}\|_{\tilde{a}}-\alpha\eta\|u_{\mathup{glo}}\|_{\tilde{a}}\|v_{\mathup{glo}}\|_{\tilde{a}}-\alpha \eta\|u_{\mathup{glo}}\|_{\tilde{a}}\|v_{\mathup{ms}}\|_{\tilde{a}}\\
&\geq\underbrace{(\beta_{\mathup{glo}}/4-9\alpha \eta/4-3\alpha\eta/2)}_{\beta_{\mathup{cem}}}\|u_{\mathup{ms}}\|_{\tilde{a}}\|v_{\mathup{ms}}\|_{\tilde{a}},
\end{align*}
where the last inequality comes from the  $ 1/2 \norm{v_{\mathup{ms}}}_{\tilde{a}} \leq \norm{v_{\mathup{glo}}}_{\tilde{a}} \leq 3/2 \norm{v_{\mathup{ms}}}_{\tilde{a}}$ and $1/2 \norm{u_{\mathup{ms}}}_{\tilde{a}} \leq \norm{u_{\mathup{glo}}}_{\tilde{a}} \leq 3/2 \norm{u_{\mathup{ms}}}_{\tilde{a}}$ by the smallness of $\eta$.
Then we obtain the Cea's lemma, if $u$ is the real solution of \cref{vf} and $u_{\mathup{ms}}$ is the multiscale solution of \cref{cemvar}, there exists a  $v_{\mathup{ms}}\in V_{\mathup{ms}}$,
\begin{align*}
\beta_{\mathup{cem}}\|u-u_{\mathup{ms}}\|_{\tilde{a}}^2\leq\mathcal{B}(u-u_{\mathup{ms}},u-u_{\mathup{ms}})&=\mathcal{B}(u-u_{\mathup{ms}},u-v_{\mathup{ms}})+\mathcal{B}(u-u_{\mathup{ms}},v_{\mathup{ms}}-u_{\mathup{ms}})\\
&\leq\alpha\|u-u_{\mathup{ms}}\|_{\tilde{a}}\|u-v_{\mathup{ms}}\|_{\tilde{a}}.
\end{align*}
\end{proof}
\begin{theorem}\label{thm:main local}
Under \cref{ass:overlapping}, there exist positive constants $\Upsilon'$ and $\Lambda^*$ such that, for any $\Upsilon \ge \Upsilon'$ and $\Lambda \ge \Lambda^*$, the multiscale solution $u_{\mathup{ms}}$ of \cref{cemvar} exists and is unique.

Moreover, for any $m \ge 1$, the following error estimate holds:
\[
\|u - u_{\mathup{ms}}\|_{\tilde{a}}
\le
C_{\mathup{cem}}(k)
\left(
\frac{1}{\varepsilon_0\sqrt{\Lambda}}
+
H^{-2}(m+1)^{d/2}\theta^{(m-1)/2}
\right)
\|f\|_{\tilde{s}^{-1}},
\]
where $C_{\mathup{cem}}(k)$ depends on $k$ and $\Upsilon^{-1}$, and $0<\theta < 1$ is a constant independent of $k$, $m$, and $H$. In particular, By 
choosing the oversampling layer $m$ such that $C_{\mathup{cem}}(m+1)^{d/2} \cdot \theta^{(m-1)/2}\cdot  \sim O(H^2)$, then the method achieves {optimal} convergence order in $H$, namely,
\[
\|u - u_{\mathup{ms}}\|_{\tilde{a}}
\approx \mathcal{O}(H).
\]
\end{theorem}
\begin{proof}
By using of Cea's Lemma in \cref{cea}, we can obtain there exists an $v_{\mathup{ms}}\in V_{\mathup{ms}}$, such that
\begin{align*}
\norm{u-u_{\mathup{ms}}}_{\tilde{a}}&\leq\frac{\alpha}{\beta_{\mathup{cem}}}\norm{u-v_{\mathup{ms}}}_{\tilde{a}}\\
&\leq \frac{\alpha}{\beta_{\mathup{cem}}}(\underbrace{\|u-u_{\mathup{glo}}\|_{\tilde{a}}}_{I_1}+\underbrace{\|u_{\mathup{glo}}-v_{\mathup{ms}}\|_{\tilde{a}}}_{I_2})\\
&\leq \frac{\alpha}{\beta_{\mathup{cem}}}\RoundBrackets*{\frac{C_1\Norm{f}_{\tilde{s}^{-1}}}{\varepsilon_0\sqrt{\Lambda}}+(c_0(\Lambda)C_{\mathup{ol}})^{1/2}C_{\mathup{po}}C_*'H^{-2}(m+1)^{d/2} \theta^{(m-1)/2}\Norm{f}_{\tilde{s}^{-1}}}\\
&\leq C_{\mathup{cem}}(k)\RoundBrackets*{\frac{\Norm{f}_{\tilde{s}^{-1}}}{\varepsilon_0\sqrt{\Lambda}}+H^{-2}(m+1)^{d/2} \theta^{(m-1)/2}\Norm{f}_{\tilde{s}^{-1}}}.
\end{align*}
For $I_1$, we utilize the global error estimate \cref{glothm}.
For $I_2$, we assume $u_{\mathup{glo}}=\mathcal{G}^{\infty}\psi$ {and} let $v_{\mathup{ms}}=\hat{\mathcal{G}}^m\psi$. After repeating the steps in \cref{cea}, we can have 
$$I_2\leq (c_0(\Lambda)C_{\mathup{ol}})^{1/2} C_*'H^{-1}(m+1)^{d/2} \theta^{(m-1)/2}\|\mathcal{G}^{\infty}\psi\|_{\tilde{a}},$$
where the estimate of $\|\mathcal{G}^{\infty}\psi\|_{\tilde{a}}$ comes from the inf-sup condition of the global problem \cref{eq:global solution}
\[
C_{\mathup{po}}H^{-1}\|w_{\mathup{glo}}\|_{\tilde{a}}\|f\|_{\tilde{s}^{-1}}\geq\|f\|_{\tilde{s}^{-1}}\|w_\mathup{glo}\|_{\tilde{s}}\geq(f,w_{\mathup{glo}})=\mathcal{B}(u_\mathup{glo},w_{\mathup{glo}})\geq\beta_{\mathup{glo}}\|\mathcal{G}^{\infty}\psi\|_{\tilde{a}}\|w_{\mathup{glo}}\|_{\tilde{a}}.
\]
\end{proof}
\section{Numerical experiments}\label{sec:Numerical}
We conduct numerical experiments on a square domain $\Omega = (0,1) \times (0,1)$. All numerical experiments are simulated by using the Python libraries NumPy and SciPy. The coefficient profile $\sigma$ and $c$ are represented by a $400 \times 400$ matrix/image, with each element corresponding to a constant value on a fine element  $\mathcal{K}_h$. The reference solution, auxiliary spaces, and multiscale bases are all calculated on the fine mesh $\mathcal{K}_h$ using the $Q_1$ finite element method with fine mesh $400\times 400$.
We measure the convergence of the proposed method using both the relative energy error and the relative $L^2$ error which are defined:
\[
  \frac{\norm{e_h}_{\tilde{a}}}{\norm{u_h}_{\tilde{a}}} \text{  and  } \frac{\norm{e_h}_{L^2(\Omega)}}{\norm{u_h}_{L^2(\Omega)}},
\]
where $u_h$ is the reference solution calculated by the $Q_1$ FEM on $\mathcal{K}_h$ or the nodal interpolation of the exact solution (if available), and $e_h$ is the error between the reference solution and the numerical solution. For the coarse mesh sizes $H$, we consider coarse meshes to be $\mathcal{K}_H$:  {$10\times 10$}, $20 \times 20$, $40 \times 40$, and $80 \times 80$.
For simplicity, we take $\mu$ as
\[
  \mu|_{K_i}=\mu_\mathup{msh}\Diam(K_i)^{-2}c|_{K_i}=24H^{-2}c|_{K_i}
\]
for all numerical experiments, as suggested in \cite{Ye2023c}.
\subsection{Flat interface model}
\begin{figure}[!ht]
\centering
\includegraphics[width=\linewidth]{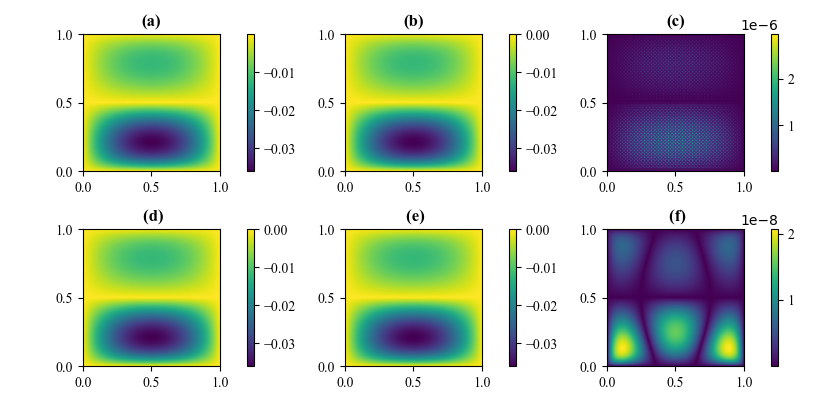} 
\caption{\textbf{(a)} The solution approximated by the CEM-GMsFEM when we choose $H=1/40$ and $m=3$.
\textbf{(d)} The solution approximated by the $Q1$ FEM when we choose $H=1/40$. {\textbf{(b)/(e)}} The exact solution. {\textbf{(c)}} The difference of the solution approximated by CEM-GMsFEM and the exact solution. {\textbf{(f)}} The difference of solution approximated by $Q_1$ FEM and exact solution.}
\label{fig:flat-solution}
\end{figure}

We initially examine a flat interface model, where the coefficient variables in $\eqref{eq:ell1}$ are defined as follows:
\begin{equation}\label{sigma}
\sigma=c=\begin{cases}
1, & y \geq 0.5\\
-3, &\text{otherwise}
\end{cases}.
\end{equation}
Then we denote $\Omega^{+}=(0,1)\times(0.5,1)$ and $\Omega^{-}=(0,1)\times(0,0.5)$. We set fixed values $(\sigma_*^{+},-\sigma_*^{-})$ to $(\sigma^+,\sigma^-)$, where both $\sigma_*^{+}$ and $\sigma_*^{-}$ are positive values. 
In the subsequent simulation, we choose $\gamma=0.5, \sigma_*^-=3$, and $\sigma_*^+=1$ from \eqref{sigma} to design the exact solution $u$ as 
\[
  u(x_1, x_2) = \begin{cases}
    -\sigma^-_*x_1(x_1-1)x_2(x_2-1)(x_2-\gamma), & \text{ in } \Omega^+, \\
    \sigma^+_*x_1(x_1-1)x_2(x_2-1)(x_2-\gamma),  & \text{ in } \Omega^-,
  \end{cases}
\]
which corresponds to a smooth source term $f$ given by
$$
f(x_1, x_2) = 
  \sigma^-_*\sigma^+_*\Big(2x_2(x_2 - 1)(x_2 - \gamma)+ x_1(x_1 - 1)(6x_2 - 2(\gamma + 1))-k^2x_1(x_1-1)x_2(x_2-1)(x_2-\gamma)\Big).
$$ 
{The theory in \cite{Chesnel2013} shows that the problem is well-posed whenever $\Upsilon \neq 1$.} For the proposed multiscale method, we fix the wavenumber $k = 4$. 
We set $l_* = 3$, meaning that the first three local eigenfunctions are computed and three multiscale basis functions are constructed for each coarse element. 
Meanwhile, the number of oversampling layers $m$ is varied from $1$ to $4$.
\begin{table}[!ht]
\caption{
The relative errors in the energy norm (columns labelled with $\norm{\cdot}_{\tilde{a}}$) 
and in the $L^2$ norm (columns labelled with $\norm{\cdot}_{L^2(\Omega)}$).
}
\label{tab:flat-error}
\centering
\resizebox{\textwidth}{!}{
\makegapedcells
\footnotesize{
\begin{tabular}{ c c c c c c c c c}
\toprule
\multirow{3}{*}{$H$} 
& \multicolumn{2}{c}{$Q_1$}  
& \multicolumn{2}{c}{$m=2$}    
& \multicolumn{2}{c}{$m=3$}  
& \multicolumn{2}{c}{$m=4$}\\
\cmidrule{2-9}
& $\norm{\cdot}_{\tilde{a}}$ & $\norm{\cdot}_{L^2(\Omega)}$
& $\norm{\cdot}_{\tilde{a}}$ & $\norm{\cdot}_{L^2(\Omega)}$
& $\norm{\cdot}_{\tilde{a}}$ & $\norm{\cdot}_{L^2(\Omega)}$
& $\norm{\cdot}_{\tilde{a}}$ & $\norm{\cdot}_{L^2(\Omega)}$ \\
\midrule            

{$\frac{1}{10}$} 
& {\num{1.258e-03}} 
& {\num{7.601e-04}}
& {\num{1.057e-01}} 
& {\num{7.661e-03}}
& {\num{4.571e-03}} 
& {\num{5.663e-05}}
& {\num{3.571e-04}} 
& {\num{3.663e-06}} \\

$\frac{1}{20}$ 
& \num{3.393e-04} & \num{1.900e-04}
& \num{6.150e-02} & \num{5.581e-03}
& \num{2.604e-03} & \num{3.220e-05}
& \num{1.082e-04} & \num{1.369e-06} \\

$\frac{1}{40}$ 
& \num{8.582e-05} & \num{4.801e-05}
& \num{2.839e-02} & \num{1.199e-03}
& \num{1.140e-03} & \num{2.652e-05}
& \num{5.660e-05} & \num{1.246e-06} \\

$\frac{1}{80}$ 
& \num{2.152e-05} & \num{1.203e-05}
& \num{1.265e-02} & \num{6.209e-04}
& \num{4.980e-04} & \num{2.206e-05}
& \num{1.713e-04} & \num{5.386e-06} \\

\bottomrule
\end{tabular}
}
}
\end{table}
\begin{figure}[!ht]
\centering
\includegraphics[width=\linewidth]{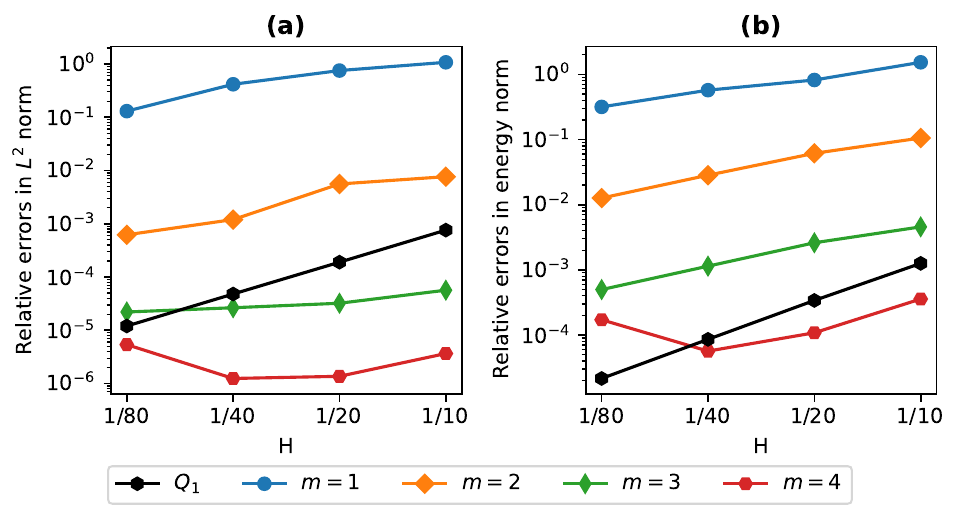}
\caption{The relative errors of the proposed method with different numbers of oversampling layers $m$ and the $Q_1$ FEM are calculated w.r.t.\ the coarse mesh size $H$.}
\label{fig:flat-error}
\end{figure}
\circled{1} From subplots \SubplotTag{(a)} to \SubplotTag{(b)} in \cref{fig:flat-error}, we observe that the convergence of the $Q_1$ FEM exhibits a linear behavior with respect to $H$ on the logarithmic scale. {This is because the interface is fully resolved by each coarse mesh and is consistent with the theoretical expectation.} \circled{2} In \cref{tab:flat-error}, we report the errors obtained under successive refinements of the coarse mesh size together with appropriately chosen oversampling layers. The first two columns correspond to the relative error in the energy norm and the relative error in the $L^2$ norm for the $Q_1$ FEM. It can be seen that for $m=3,4$, the proposed multiscale method achieves higher accuracy than the $Q_1$ FEM. For instance, when $H = 1/20$, the relative energy error of the multiscale method is $0.0108\%$, whereas the $Q_1$ FEM yields $0.0339\%$. \circled{3} \Cref{fig:flat-error} further illustrates the errors of the CEM-GMsFEM for different coarse mesh sizes and oversampling layers. However, for a fixed $m$, the error does not always decrease monotonically with respect to $H$, as shown in subplots \SubplotTag{(a)} and \SubplotTag{(b)}. This behavior differs from that of classical FEMs but is typical for CEM-GMsFEM. The reason is that a larger coarse mesh size $H$ leads to a larger oversampling region, which captures more local information, as indicated in \cref{thm:main local}. \circled{4} In addition, we visualize the exact solution in \cref{fig:flat-solution}-{\textbf{(b)/(e)}}, the multiscale solutions in \cref{fig:flat-solution}{\textbf{(a)}}, and the $Q_1$ FEM solution in \cref{fig:flat-solution}-{\textbf{(d)}}. These comparisons further demonstrate the efficiency of the proposed method.
\subsection{Random inclusion model}
\begin{figure}
\centering
\includegraphics[width=\linewidth]{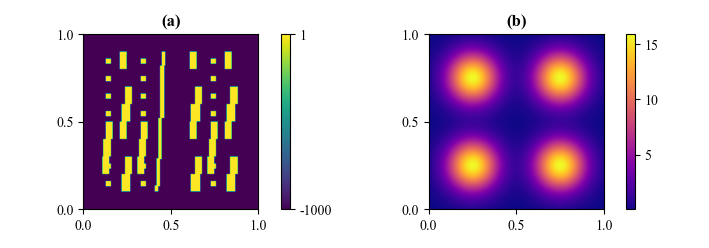}
\caption{\textbf{(a)}  Coefficients with $(\sigma_*^+,\sigma_*^-)=(1,10^3)$; \textbf{(b)} Source term based on Gaussian functions.}
\label{fig:High-contrast}
\end{figure}
\begin{figure}[!ht]
\centering
\includegraphics[width=\linewidth]{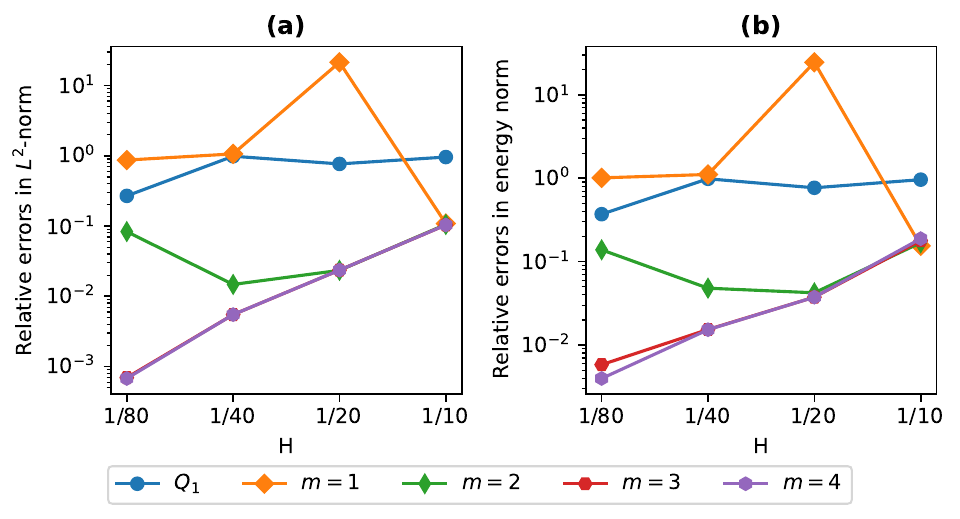}
\caption{The relative errors of the proposed method with different numbers of oversampling layers $m$ and the $Q_1$ FEM are calculated w.r.t.\ the coarse mesh size $H$.}
\label{fig:ex2derror}
\end{figure}
\begin{figure}[!ht]
\centering
\includegraphics[width=\linewidth]{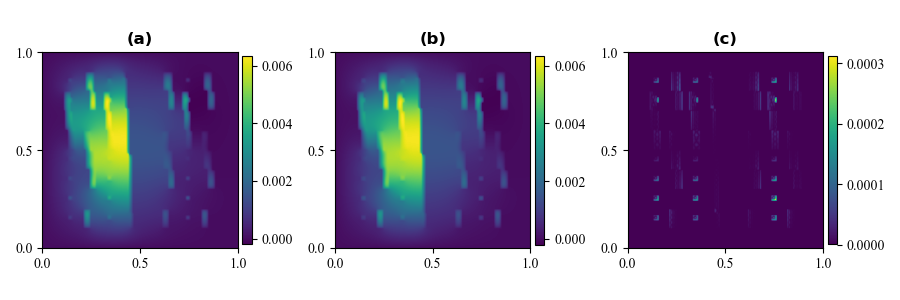}
\caption{\textbf{(a)} The solution approximated by the CEM-GMsFEM when we choose $H=1/40$ and $m=3$.
\textbf{(b)} The solution approximated by the $Q1$ FEM when we choose $H=1/40$.}
\label{fig:ex20}
\end{figure}
\begin{figure}[ht]
    \centering
    \includegraphics[width=\linewidth]{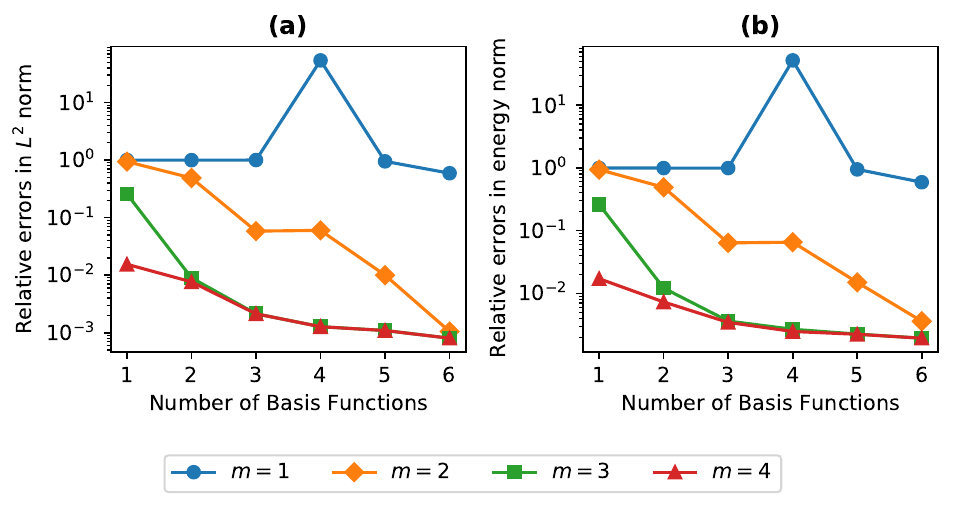}
    \caption{The relative errors of the proposed method with different numbers of basis functions when $m=3$ and $H=1/20$.}
    \label{fig:Ex2basis}
\end{figure}

In this subsection, we consider a random inclusion model that is utilized in several multiscale methods as a showcase of the capability of handling nonperiodic coefficient profiles \cite{Zhao2020,Poveda2024}.
The subdomains $\Omega^+$ and $\Omega^-$ are demonstrated in \cref{fig:High-contrast}-\SubplotTag{(a)}, and $\sigma$ is determined again by $(\sigma^+_*, \sigma^-_*)$.
We consider $(\sigma^+_*, \sigma^-_*)=(1, 10^{3})$ and 
set the source term $f$ to be the following Gaussian function 
$$f(x_1,x_2)=\frac{1}{\sigma\sqrt{2\pi}}\exp(-r^2/2\sigma^2)$$
where $r$ is the distance from the point $(x_1,x_2)$ to the center and $\sigma$ is the standard deviation. The source function is depicted in \cref{fig:High-contrast}-\SubplotTag{(b)}.  {Although $\Upsilon = 10^{-3}$ is small in this case, as noted in Remark \ref{remark1}, one may equivalently define the reversed flipping operator so that the \texttt{T}-coercivity argument remains valid when $\Upsilon$ is sufficiently small, provided that $\Upsilon \neq 1$.} For the proposed multiscale method, we fix $k = 4$ and $l_* = 3$. \circled{1} In \cref{fig:High-contrast}, the standard $Q_1$ FEM fails to provide satisfactory accuracy in both the relative $L^2$ norm and the energy norm. In contrast, for the multiscale method, when $m=4$, clear and stable convergence is observed. The classical CEM-GMsFEM \cite{Jin2025, Ye2023c} is proven to be effective in handling long and high-contrast channels, and the proposed method inherits this advantage. \circled{2} In \cref{tab:ex2}, we report the errors obtained under successive refinements of the coarse mesh size with different choices of oversampling layers. As the coarse mesh size $H$ decreases, both the relative $L^2$ error and the energy error of the multiscale solution exhibit clear convergence toward the fine-grid reference solution computed by the $Q_1$ FEM. \circled{3} As shown in \cref{fig:ex2derror}, an interesting phenomenon can be observed: for $H = 1/20$ and $H = 1/40$, the error does not decrease significantly when increasing $m$ from $3$ to $4$. This indicates that the overall error is primarily limited by the coarse mesh size $H$, and that $m=3$ is already sufficient to achieve an accurate approximation. To further illustrate this observation, we compare the multiscale solution and the $Q_1$ FEM solution for $H = 1/40$ and $m=3$ in \cref{fig:ex20}. The two solutions are visually almost indistinguishable, with a difference of only $1.60\%$. \circled{4} We further investigate the influence of the number of basis functions on the relative $L^2$ and energy errors for different oversampling layers in the case $H = 1/20$. The results, presented in \cref{fig:Ex2basis}, show that increasing the number of basis functions leads to a noticeable reduction in both error norms. In particular, a clear turning point appears between $l_* = 3$ and $l_* = 4$ when $m = 1$. Considering the balance between computational cost and desired accuracy, we therefore select $l_* = 3$. These results highlight the importance of choosing an appropriate number of basis functions to enhance the accuracy of multiscale methods.
\begin{table}[!ht]
\caption{The relative errors in the energy norm (in the columns labelled with $\norm{\cdot}_{\tilde{a}}$) and in the $L^2$ norm (in the columns labelled with $\norm{\cdot}_{L^2(\Omega)}$).
  }\label{tab:ex2}
\centering
\resizebox{\textwidth}{!}{
\makegapedcells
\footnotesize{
\begin{tabular}{ c c c c c c c c c}
\toprule
\multirow{3}{*}{$H$} & \multicolumn{2}{c}{$m=1$}  & \multicolumn{2}{c}{$m=2$}    & \multicolumn{2}{c}{$m=3$}  & \multicolumn{2}{c}{$m=4$}\\
\cmidrule{2-9}
 & $\norm{\cdot}_{\tilde{a}}$ & $\norm{\cdot}_{L^2(\Omega)}$ & $\norm{\cdot}_{\tilde{a}}$ & $\norm{\cdot}_{L^2(\Omega)}$ & $\norm{\cdot}_{\tilde{a}}$ & $\norm{\cdot}_{L^2(\Omega)}$ & $\norm{\cdot}_{\tilde{a}}$ & $\norm{\cdot}_{L^2(\Omega)}$ \\
\midrule            
{$\frac{1}{10}$}
& {\num{1.547e-01}} 
& {\num{1.084e-01}}
& {\num{1.686e-01}} 
& {\num{1.053e-01}}
& {\num{1.786e-01}} 
& {\num{1.035e-01}}
& {\num{1.886e-01}} 
& {\num{1.035e-01}} \\

$\frac{1}{20}$
& \num{2.457e+01} & \num{2.135e+01}
& \num{4.231e-02} & \num{2.340e-02}
& \num{3.750e-02} & \num{2.357e-02}
& \num{3.749e-02} & \num{2.357e-02} \\

$\frac{1}{40}$
& \num{1.110e+00} & \num{1.062e+00}
& \num{4.804e-02} & \num{1.475e-02}
& \num{1.536e-02} & \num{5.488e-03}
& \num{1.526e-02} & \num{5.485e-03} \\

$\frac{1}{80}$
& \num{1.011e+00} & \num{8.677e-01}
& \num{1.385e-01} & \num{8.327e-02}
& \num{5.810e-03} & \num{6.962e-04}
& \num{3.971e-03} & \num{6.771e-04} \\

\bottomrule
\end{tabular}
}
}
\end{table}

\subsection{2D negative-index materials}
\begin{figure}[!ht]
\centering
\includegraphics[width=\linewidth]{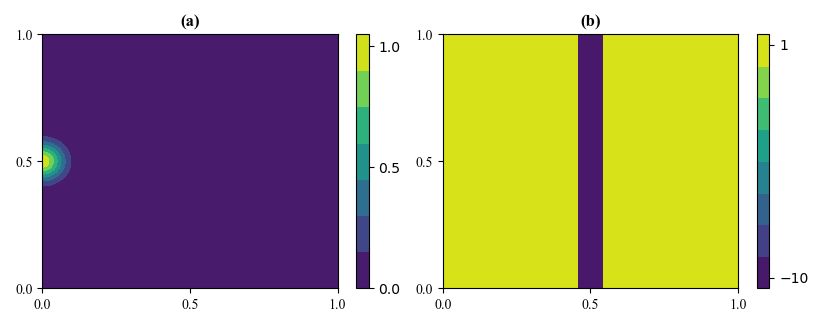}
\caption{\textbf{(a)} Gaussian Beam Source centered at $(0,0.5)$.
\textbf{(b)} The Negative index Materials.}
\label{fig:2dnim}
\end{figure}

\begin{figure}[hbtp]
\centering
\includegraphics[width=\linewidth]{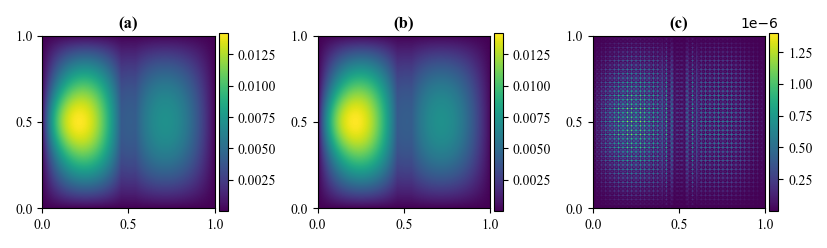}
\caption{\textbf{(a)} The solution approximated by the CEM-GMsFEM when we choose $H=1/40$ and $m=3$. \textbf{(b)} The solution approximated by the $Q1$ FEM when we choose $H=1/40$.}
\label{fig:2dnim11}
\end{figure}

\begin{figure}[hbtp]
\centering
\includegraphics[width=\linewidth]{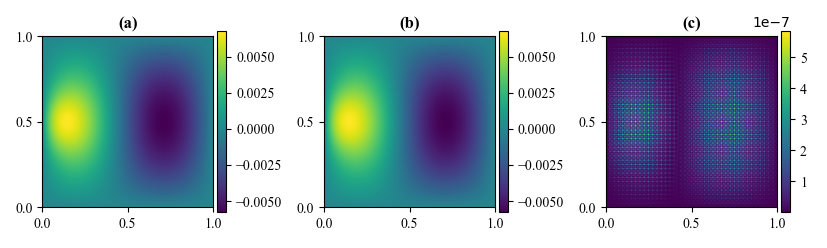}
\caption{\textbf{(a)} The solution approximated by the CEM-GMsFEM when we choose $H=1/40$ and $m=3$.
\textbf{(b)} The solution approximated by the $Q1$ FEM when we choose $H=1/40$.}
\label{fig:2dnim10}
\end{figure}

In this section, we provide an example demonstrating the interactions of Gaussian beams with 2D negative-index materials \cite{Li2008, ziolkowski2003}. The coefficients in $\eqref{eq:ell1}$ are defined as 
\begin{equation*}
\sigma=c=
    \begin{cases}
1,&\text{if}\,\,  x<11/24\,\,\text{or}\,\, x>13/24\\
-10,&\text{otherwise}
    \end{cases},
\end{equation*}
and are depicted in \cref{fig:2dnim}-{\SubplotTag{(b)}.} Consequently, the subdomains $\Omega^+$ and $\Omega^-$ are demonstrated in \cref{fig:2dnim}-{\SubplotTag{(b)},}  and $\sigma$ is determined again by $(\sigma^+_*, \sigma^-_*)=(1,10)$. The source wave is set as 
$$f(x,y)=\exp(-((x-x_0)^2+(y-y_0)^2)/(2\text{waist}^2)),$$
where beam center $(x_0,y_0)=(0,0.5)$ and  $\text{waist}=0.05$, which is depicted in \cref{fig:2dnim}-{\SubplotTag{(a)}.} 
\begin{table}[!ht]
\caption{
The relative errors in the energy norm (in the columns labelled with $\norm{\cdot}_{\tilde{a}}$) and in the $L^2$ norm (in the columns labelled with $\norm{\cdot}_{L^2(\Omega)}$).
  }\label{tab:2d-nim}
\centering
\resizebox{\textwidth}{!}{
\makegapedcells
\footnotesize{
\begin{tabular}{ c c c c c c c c c}
\toprule
\multirow{3}{*}{$H$} & \multicolumn{2}{c}{$m=1$}  & \multicolumn{2}{c}{$m=2$}    & \multicolumn{2}{c}{$m=3$}  & \multicolumn{2}{c}{$m=4$}\\
\cmidrule{2-9}
 & $\norm{\cdot}_{\tilde{a}}$ & $\norm{\cdot}_{L^2(\Omega)}$ & $\norm{\cdot}_{\tilde{a}}$ & $\norm{\cdot}_{L^2(\Omega)}$ & $\norm{\cdot}_{\tilde{a}}$ & $\norm{\cdot}_{L^2(\Omega)}$ & $\norm{\cdot}_{\tilde{a}}$ & $\norm{\cdot}_{L^2(\Omega)}$ \\
\midrule            
{$\frac{1}{10}$}
& {\num{1.381e-01}} 
& {\num{8.693e-02}}
& {\num{5.410e-02}} 
& {\num{7.591e-03}}
& {\num{5.382e-02}} 
& {\num{7.565e-03}}
& {\num{5.382e-02}} 
& {\num{7.565e-03}} \\

$\frac{1}{20}$ 
& \num{2.532e-01} & \num{1.900e-01}
& \num{1.204e-02} & \num{7.733e-04}
& \num{6.251e-03} & \num{3.389e-04}
& \num{6.231e-03} & \num{3.374e-04} \\

$\frac{1}{40}$ 
& \num{6.290e-01} & \num{3.300e-01}
& \num{2.312e-02} & \num{1.724e-03}
& \num{1.314e-03} & \num{3.613e-05}
& \num{9.143e-04} & \num{2.438e-05} \\

$\frac{1}{80}$ 
& \num{6.504e-01} & \num{8.072e-01}
& \num{5.198e-02} & \num{9.251e-03}
& \num{2.390e-03} & \num{4.200e-05}
& \num{1.676e-04} & \num{2.276e-06} \\
\bottomrule
\end{tabular}
}
}
\end{table}
\begin{figure}[!ht]
\centering
\includegraphics[width=\linewidth]{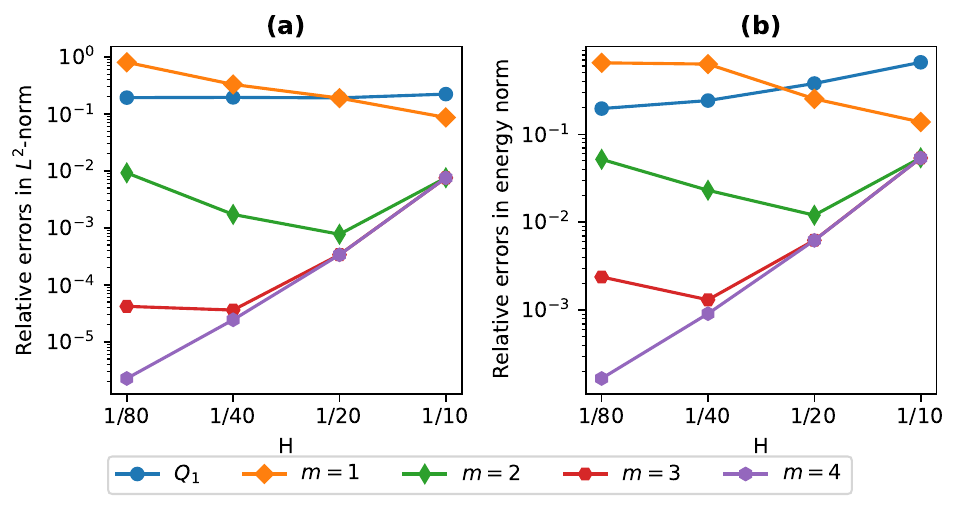}
\caption{The relative errors of the proposed method with different numbers of oversampling layers $m$ and the $Q_1$ FEM are calculated w.r.t.\ the coarse mesh size $H$.}
\label{fig:2dnimerror}
\end{figure} We proceed by presenting the numerical errors of the proposed method under the same parameter settings as in the previous example. \circled{1} In \cref{tab:2d-nim} and \cref{fig:2dnimerror}, we compare the errors for different coarse mesh sizes and oversampling layers. Since an analytical solution is not available, the fine-grid $Q_1$ FEM solution is used as the reference solution. We observe a behavior similar to that in the second example: the standard $Q_1$ FEM does not exhibit the desired convergence, whereas the proposed multiscale method demonstrates improved convergence rates and higher accuracy. \circled{2} To illustrate the exotic properties of the materials under consideration, the phenomenon of negative refraction is depicted in \cref{fig:2dnim11}. In this figure, the intensity is represented by a color gradient, where high values are shown in yellow and low values in dark blue. As the beam propagates through the slab region $[1/3, 2/3]$, a vacuum-like region can be observed in \cref{fig:2dnim11}-\SubplotTag{(a)} , clearly demonstrating the refraction characteristics of the negative-index metamaterial (NIM) slab. This effect is further illustrated in \cref{fig:2dnim10}, where we plot the solution in the non-sign-changing domain for $H = 1/40$ and $m = 3$. \circled{3} This example can be regarded as an in-silico simulation of wave propagation in metamaterials. The behavior of the wave as it passes through the negative-index medium further demonstrates the effectiveness and robustness of the proposed multiscale method.
\section{Conclusion}\label{sec:conclusions}
We have proposed a multiscale computational method CEM-GMsFEM, for solving electromagnetic wave propagation problems. Due to the sign-changing coefficients in the model problem, the coercivity of the bilinear form cannot be guaranteed. 
To address this issue, we employ the \texttt{T}-coercivity framework, replacing the coefficient with its absolute value and demonstrating that this modification is consistent with the \texttt{T}-coercivity theory as well as the required resolution condition. Under several technical assumptions, we establish rigorous error estimates and obtain optimal convergence rates.  The numerical experiments confirm the effectiveness and robustness of the proposed method. For the flat interface model, the multiscale method exhibits greater flexibility and higher accuracy than the standard $Q_1$ finite element method when appropriate oversampling layers are employed. In the second example, involving a random inclusion model with high-contrast coefficients, the proposed method again demonstrates superior convergence compared to the $Q_1$ FEM. Moreover, we observe that using a small number of oversampling layers can still provide accurate approximations while significantly reducing computational cost. 
We also discuss practical strategies for selecting an appropriate number of multiscale basis functions to achieve optimal performance. In the final example, we simulate wave propagation through metamaterials and observe the formation of a vacuum region, illustrating the remarkable ability of metamaterials to manipulate wave behaviour.

Multiscale basis functions are widely recognized as effective coarse spaces in the construction of multilevel domain decomposition or multigrid  preconditioners \cite{Ye2024a, Ye2025}. In this framework, the favorable properties of multiscale computational methods, such as robustness with respect to high contrast, can be 
leveraged to accelerate iterative solvers. In the context of the present problem, three related research directions have emerged: the sign-changing Laplace-type PDE, the Helmholtz-type PDE, and the {curl--curl-type PDE}. The mathematical theory for the first class 
is relatively well established, while the latter two—particularly { curl--curl case—remain less understood and are still under active development \cite{Rihan2025}, especially for the 3D simulations.} We emphasize that, due to the non-coercive nature of these problems, 
the direct application of conventional multigrid preconditioners is not straightforward and may even fail. One possible strategy is to reformulate the original problem into a coercive one by enriching the 
approximation space. In such a framework, multiscale basis functions can be incorporated as a key component of the overall solver. A satisfactory resolution of this issue would substantially enlarge the feasible computational scale of current NIM simulations.

\end{document}